\documentclass[a4paper,10pt]{article}
\usepackage{amssymb,amsmath,epsfig,theorem}
\usepackage[english]{babel}

\input xy
\xyoption{all}

\newtheorem{prop}{Proposition}[section]
\newtheorem{theor}[prop]{Theorem}
\newtheorem{lemm}[prop]{Lemma}
\newtheorem{defin}[prop]{Definition}
\newtheorem{corollary}[prop]{Corollary}

\begin{document}

\title{On the analytical invariance of the semigroups of a quasi-ordinary 
hypersurface singularity}
\author{Patrick Popescu-Pampu}
\date{}

\maketitle
\thispagestyle{empty}
%%%%%%%%%%%%%%%%%%%%%%%%%%%%%%%%%%%%%%%%%%%%%%%%%%%%%%%%%%%%
%%%  Abstract  %%%
%%%%%%%%%%%%%%%%%%
\begin{abstract}{
We associate to any irreducible germ $\mathcal{S}$ 
of complex quasi-ordinary hypersurface an analytically invariant
semigroup.  
We deduce a direct proof (without passing through their 
embedded topological invariance) 
of the analytical invariance of the normalized 
characteristic exponents. These exponents generalize  
the generic Newton-Puiseux exponents of plane curves. Incidentally,
we give a toric description of the normalization morphism of the germ
$\mathcal{S}$.} 
\end{abstract}

{\small 2000 \textit{Mathematics Subject Classification.} Primary
    32S10; Secondary 14M25.}

\par\medskip\centerline{\rule{2cm}{0.2mm}}\medskip
\setcounter{section}{0}
%%%%%%%%%%%%%%%%%%%%%%%%%%%%%%%%%%%%%%%%%%%%%%%%%%%%%%%%%%%%
%%%  Main text (in English)  %%%
%%%%%%%%%%%%%%%%%%%%%%%%%%%%%%%%
%\setcounter{section}{1}

\section{Introduction}

In this paper we generalize to arbitrary dimensions results obtained
first for surfaces in \cite{PP 01} and published in \cite{PP 02}.  

A classical way to study an irreducible germ $C$ of complex analytical 
plane curve 
is to introduce its Newton-Puiseux series in some coordinate system 
$(X,Y)$, which allows to define its so-called \textit{characteristic
  (Newton-Puiseux) exponents}. 
If the coordinates are \textit{generic} - which means that the
$Y$-axis and the 
embedded reduced tangent cone of the curve are transversal - the 
characteristic exponents 
do not depend on them and their collection is a complete invariant of 
the embedded topological type of the curve.

Another way to study the germ $C$ is to associate to it a \textit{semigroup}
$\Gamma(C)$. We recall two ways of doing this. Both of them yield 
isomorphic (abstract) semigroups:

1) Take the values by the canonical valuation of the elements of the
local algebra $\mathcal{A}$ of $C$. In other words, take the orders of
vanishing of $\nu^{*}(h)$ at the base point of  
$\overline{C}$, where $\nu: \overline{C} \rightarrow C$ is the normalization 
morphism of $C$ and $h$ varies through $\mathcal{A}$. They form a
sub-semigroup of $(\mathbf{N},+)$ which obviously   
depends only on the analytical type of the germ $C$.

2) Take the orders of the series 
$h(\xi)$, where $h$ varies in $\mathbf{C}\{X,Y\}$ and $\xi$ is a 
Newton-Puiseux series of $C$ in the coordinates $(X,Y)$. They form a
sub-semigroup of $(\mathbf{Q}_{+},+)$ which can be expressed in terms
of the characteristic exponents of $C$. 

Seen as an abstract 
semigroup, $\Gamma(C)$ 
is also a complete invariant of the embedded topological type of $C$. 
For the preceding claims, see \cite{Z 86} and \cite{PP 99}. In
\cite{PP 99}, we noticed that from the isomorphism of the two
semigroups one can deduce the analytical invariance of the (generic) 
characteristic exponents of $C$. 
\medskip

Here we extend this
idea to a class of higher-dimensional hypersurface germs, the so-called
\textit{quasi-ordinary} ones, for which a generalization of the
characteristic exponents can be defined. If $\mathcal{S}$ is an
irreducible germ of hypersurface of dimension $d$, such exponents are
associated to any 
\textit{quasi-ordinary projection}:
$$\psi:\mathcal{S}\rightarrow \mathbf{C}^{d},$$
which is by definition a finite morphism whose discriminant locus
is contained in a hypersurface with normal crossings.

Lipman generalized the notion of generic Newton-Puiseux exponents of
plane curves, by
defining the \textit{normalized} characteristic exponents of
$\mathcal{S}$ (see section \ref{basdef}).  The irreducible
quasi-ordinary germ of hypersurface $\mathcal{S}$ 
being given, there is always a quasi-ordinary projection which has
moreover normalized characteristic exponents.   

It is a natural question to study the degree of invariance (analytic
or topological) of the
normalized characteristic exponents. In \cite{G 86}, Gau proves their
embedded topological invariance when $\mathcal{S}$ is a germ of
surface. Then, Gau and Lipman \cite{G 88}, \cite{L 88} generalize
this result to arbitrary dimensions. 

In \cite{GP 00}, \cite{GP 02}, 
Gonz{\'a}lez P{\'e}rez  generalizes the second of the 
constructions presented above of the semigroup of a plane curve 
to the case of an irreducible quasi-ordinary hypersurface germ
$\mathcal{S}$. He starts from a fixed quasi-ordinary polynomial
$f$ which defines $\mathcal{S}$ (see definition \ref{qopol}).  
Instead of the order of a series in one variable, he uses the set of
vertices of the Newton polyhedron of a series in various variables. He
defines the semigroup $\Gamma(f)$ as the set of vertices of the Newton
polyhedra of the fractional series $h(\xi)$, where $h$ varies in
$\mathbf{C}\{X_{1},...,X_{d}\}[Y]-(f)$ and $\xi$ is a fractional series
representing a root of $f$ (see section \ref{basdef}). 
The same semigroup is obtained
if one considers only functions $h$ such that $h(\xi)$ has a
dominating exponent (see section \ref{extrsg}).  
Using the embedded topological invariance of 
the normalized characteristic exponents, he shows that up to isomorphism, this 
semigroup does not depend on the quasi-ordinary polynomial $f$. 
Moreover, it is 
a complete invariant of the embedded topological type of the germ, 
and a fortiori it is an \textit{analytical invariant} of the germ. 

\medskip

Here we avoid using Lipman-Gau's results. Instead, we generalize the
first definition given in the introduction.

If $f$ is a quasi-ordinary polynomial defining $\mathcal{S}$, 
we remark that one can take as the root $\xi$ of $f$ the restriction
$Y\mid_{\mathcal{S}}$, as $f(Y) \mid_{\mathcal{S}}=0$. With this
choice of root, we have the equality: 
\begin{equation} \label{relfund}
\psi^{*}(h(\xi))=h\mid_{\mathcal{S}}.
\end{equation} 
This remark
is the starting point of our method of construction of an intrinsic
semigroup using the elements of the algebra $\mathcal{A}$.

In dimension
$d \geq 2$, there is no canonical valuation associated to
$\mathcal{S}$. We generalize definition 1 of the introduction by 
constructing a canonical morphism $\theta:(\overline{\mathcal{R}},P)\rightarrow
(\mathcal{S},0)$ whose source is a smooth space, and a canonical
divisor with normal crossings $\overline{\mathcal{H}}$ on
$\overline{\mathcal{R}}$ at $P$. Then we restrict to those functions
$h \in \mathcal{A}$ such that the 
components of the divisor $(\theta ^{*}(h))$ are either components of
$\overline{\mathcal{H}}$ or do not contain the intersection of its 
components (we say then that $\theta ^{*}(h)$ \textit{has a dominating
exponent with respect to $\overline{\mathcal{H}}$ at $P$}, this
exponent being the tuple
formed by the orders of vanishing of $\theta ^{*}(h)$ along the
components of $\overline{\mathcal{H}}$). 
The set of these dominating exponents obviously forms a semigroup of rank 
not greater
than the number of components of $\overline{\mathcal{H}}$ at $P$. We
denote it by $\Gamma_{P}'(\mathcal{S})$.

To construct the morphism $\theta$, the idea is to use the structure
of the couple $(\mathcal{S}, \mbox{Sing}(\mathcal{S}))$. For $d\geq
2$, one cannot take simply as morphism $\theta$ a 
normalization  $\nu: \overline{\mathcal{S}}\rightarrow \mathcal{S}$,
as $\overline{\mathcal{S}}$ is not necessarily smooth. So, in order to
get a smooth source, we compose $\nu$ with a finite morphism $\mu:
\mathcal{R}\rightarrow \overline{\mathcal{S}}$ (see section 11), that
we call an \textit{orbifold mapping}. We construct the hypersurface
with normal crossings by looking at the preimage $(\nu\circ
\mu)^{-1}(\mathrm{Sing}(\mathcal{S}))$. This preimage is not necessarily a
hypersurface, as $\mathrm{Sing}(\mathcal{S})$ may have components of
codimension 2 (see section 7). Therefore, we look only at its components of
codimension 1 in $\mathcal{S}$. Let $s$ be their number.

If $s=d$, we take $\theta:= \nu \circ \mu, \: P:= \theta ^{-1}(0)$ and
$\overline{\mathcal{H}}:=  \theta
^{-1}(\mathrm{Sing}(\mathcal{S}))$. Then $\theta
^{*}(h\mid_{\mathcal{S}})$ has a d.e. with respect to
$\overline{\mathcal{H}}$ at $P$, whenever $h(\xi)$ has a d.e. The
relation (\ref{relfund}) allows to construct a
morphism $\Phi_{P}:\Gamma(f) \rightarrow \Gamma'_{P}(\mathcal{S})$. In
this case, our main theorem says that $\Phi_{P}$ \textit{realizes an
isomorphism of semigroups} (theorem \ref{mths}).

If $s<d$, the situation is more complex. Then (see section
\ref{canoblo}),  we construct a
third morphism 
$\eta: \overline{\mathcal{R}}\rightarrow \mathcal{R}$ as a composition
of blowing-ups of smooth centers, determined canonically by the
structure of $(\mathcal{R},(\nu\circ
\mu)^{-1}(\mathrm{Sing}(\mathcal{S})))$. We do this third step in order
to get more components of codimension 1 of the preimage of
$\mbox{Sing}(\mathcal{S})$ passing through a same point.  We arrive at
germs $\overline{\mathcal{H}}_{P}$ having $c'$ components, where $c'$
denotes \textit{the reduced equisingular dimension} of $\mathcal{S}$,
defined in section \ref{redsg}.

In order to get an isomorphism between $\Gamma_{P}'(\mathcal{S})$ and
a semigroup generalizing construction 2 for plane curves, we are
obliged to modify Gonz{\'a}lez P{\'e}rez' definition in order to get first
of all an equality of ranks. We define (see section  \ref{redsg}) the
\textit{reduced semigroup} $\Gamma'(f)$ of $f$ by taking the vertices
of the Newton polyhedron of $h(\xi)$, where $h(\xi)$ is
now expressed as a fractional series in the first $c'$ variables
$X_{1},...,X_{c'}$. It is a semigroup of rank $c'$.

We define then a finite set $\overline{\mathcal{P}}\subset
\overline{\mathcal{H}}$ such that for any $P \in
\overline{\mathcal{P}}$, the situation is analogous to the one
explained before in the case $s=d$.  Namely, the function $\theta
^{*}(h\mid_{\mathcal{S}})$ has a d.e. with respect to
$\overline{\mathcal{H}}$ at $P$, whenever $h(\xi)$ has a d.e. The
relation (\ref{relfund}) allows to construct a
morphism $\Phi_{P}: \Gamma'(f)\rightarrow \Gamma_{P}'(\mathcal{S})$
(see section 13).

 Our main theorem is:
\medskip

\noindent \textbf{Theorem  \ref{isosgsg}} \textit{The morphism
  $\Phi_{P}$ realizes an isomorphism between the semigroups $\Gamma'(f)$ and 
$\Gamma'_{P}(\mathcal{S})$.}
\medskip

The important point for the proof  is that the morphisms $\nu, \mu,
\eta$ have \textit{toric representatives}. In particular, we think our
toric construction of a normalization morphism (see section
\ref{norgerm}) has independent interest. The needed background of
toric geometry is recalled in section \ref{remtor}.

As $\Gamma'(f)$ depends on $f$ but not on $P$ and
$\Gamma'_{P}(\mathcal{S})$ just the other way round, we see that
up to isomorphism, \textit{the semigroup obtained like this is an analytic
invariant of $\mathcal{S}$}. In section \ref{pfcor}, we deduce from this fact:

\medskip
\noindent \textbf{Corollary  \ref{invchar}} 
\textit{
The normalized characteristic exponents are analytical invariants of
the germ $\mathcal{S}$.} 
\medskip

 Other proofs of Corollary  \ref{invchar} were given in the case of 
surfaces by Lipman 
\cite{L 65}, \cite{L 83} and Luengo Velasco \cite{LV 79}, \cite{LV 83}. 
Their method was to compare the normalized characteristic exponents with 
a canonical resolution procedure by blowing-ups of smooth
centers. Such a canonical procedure by blowing-ups of smooth centers
in a way deducible only from the characteristic exponents is presently 
unavailable in higher dimensions. In the case of surfaces, an embedded 
resolution of this kind was obtained by Ban and McEwan in \cite{BM 00}.
\medskip

At a first reading, one should understand the proof of theorem
\ref{mths}. The needed material is presented in the sections which
precede it. In the sections \ref{canoblo}-\ref{pfcor} we present the
modifications needed to prove the general case. We end the paper by a
comparison with the case \cite{PP 02} of dimension 2.
\medskip

\textbf{Acknowledgements: } I am grateful to Bernard Teissier for his
suggestions and 
many encouragements, to Pedro D. Gonz{\'a}lez P{\'e}rez for his careful reading and
his suggestions
concerning the structure of this paper, to Evelia Garc{\'\i}a
Barroso and Olivier Piltant for their comments.

\section{The quasi-ordinary hypersurface germs\\ and their characteristic
  exponents} \label{basdef}

In this section we introduce the basic objects of study of this
article.

If $P$ is a point on a complex space $\mathcal{V}$, we denote by 
$\mathcal{O}_{\mathcal{V},P}$ the local algebra of $\mathcal{V}$ at
$P$. If $\mathcal{H}$ is a hypersurface on an analytical manifold
$\mathcal{V}$ and $P$ is a point of  $\mathcal{H}$, we say that
$\mathcal{H}$ has \textit{normal crossings} at $P$ if there are local
coordinates for $\mathcal{V}$ at $P$ such that the germ
$\mathcal{H}_{P}$ of $\mathcal{H}$ at $P$ is the
union of some hypersurfaces of coordinates. If $\phi$ is an analytic
function defined on a complex space, we 
denote by $\mathcal{Z}(\phi)$ its (reduced) zero-locus. 

In the sequel we will
denote with the same letter a germ and a sufficiently small
representative of it. It will be deduced from the context if one deals
with one or the other notion.  

If $a \in \mathbf{R}$, we denote by $[a]$ its integral part. 
If $F\in K[Y]$ is a polynomial with coefficients in a field $K$, we denote by
$d_{Y}(F)$ its degree and by $R(F)$ the set of its roots in an
algebraic closure of the basis field, taken with multiplicities. 
If $E$ is a finite set,
we denote by $\mid E \mid$ its cardinal. If $a,b \in \mathbf{Z}$ with 
$a \leq b$, we denote $[[a,b]]:= \{a, a+1,...,b\}$.

Let $d\geq 1$ be an integer. Define the 
\textit{algebra of fractional series} $\widetilde{\mathbf{C}\{ X\}}:= 
\mbox{lim}_{\stackrel{\rightarrow}{N \geq 0}}
\mathbf{C}\{ X_{1}^{\frac{1}{N}},...,X_{d}^{\frac{1}{N}}\}$, where 
$X:=(X_{1},...,X_{d})$. 
If $m=(m_{1},...,m_{d})\in \mathbf{Q}_{+}^{d}$, we denote $X^{m}:=
X_{1}^{m_{1}}\cdots X_{d}^{m_{d}}$. 
If $\eta \in \widetilde{\mathbf{C}\{ X\}}$ can be written 
$\eta=X^{m}u(X)$, with $ m \in \mathbf{Q}_{+}^{d}$ and 
$u \in \widetilde{\mathbf{C}\{ X\}}, \: u(0,...,0)\neq 0$, we say that 
$\eta$ has a \textit{dominating exponent} (abbreviated "d.e."), this
exponent being denoted by  $v_{X}(\eta):= m$.

If $\eta \in \widetilde{\mathbf{C}\{ X \}}$, we define its 
\textit{Newton polyhedron} $\mathcal{N}_{X}(\eta)$ to be the convex hull in 
$\mathbf{R}^{d}$ of the set $\mbox{Supp}_{X}(\eta)+ \mathbf{R}_{+}^{d}$, 
where $\mbox{Supp}_{X}(\eta)$ denotes the support of $\eta$ written 
as a series in the variables $X$. If $\eta$ has a d.e., then
$\mathcal{N}_{X}(\eta)=\{v_{X}(\eta)\}+\mathbf{R}_{+}^{d}$, which
shows that the Newton polyhedron is a generalization of the dominating
exponent. 

\begin{defin} \label{qobj}
Let $\mathcal{A}$ be a reduced equi-dimensional local complex-analytical 
algebra of dimension $d$ and $(\mathcal{S},0)$ a germ of complex space 
such that 
$\mathcal{O}_{\mathcal{S},0}\simeq \mathcal{A}$. 
The algebra $\mathcal{A}$ and the germ $(\mathcal{S},0)$ are called 
\textbf{quasi-ordinary} if there exists a finite morphism $\psi$ from 
$(\mathcal{S},0)$ to a smooth space of the same dimension, 
whose discriminant locus is contained in a hypersurface with 
normal crossings. Such a morphism $\psi$ is also called 
\textbf{quasi-ordinary}.
\end{defin}

All germs of curves are quasi-ordinary with respect to any finite
morphism whose target is a smooth curve. 

Quasi-ordinary germs appear 
naturally in the Jung method of resolution of the singularities of 
a germ by embedded resolution of the discriminant locus of a finite 
morphism from the germ to a smooth space of same dimension 
(see the original article Jung \cite{J 08} and a more recent
presentation by Laufer \cite{L 71} for the case of dimension
2). Zariski \cite{Z 67} gave an alternative method of resolution of
singularities of surfaces that need a study of quasi-ordinary germs. 
Quasi-ordinary 
hypersurface germs were first systematically studied by Lipman \cite{L 65} 
when $d=2$, see also the survey \cite{L 83}. This study was extended to 
any $d \geq 2$ in Lipman \cite{L 88}.

In the special case in which $\mathcal{A}$ is of embedding dimension $d+1$, 
one can find local coordinates $X$ on the target space of $\psi$ 
such that the discriminant locus of $\psi$ is contained in
$\mathcal{Z}(X_{1}\cdots X_{d})$ and an element $Y$ in the maximal ideal 
of $\mathcal{A}$ such that $(\psi, Y)$ embeds $(\mathcal{S},0)$ in 
$\mathbf{C}^{d}\times \mathbf{C}$. So $\psi$ appears as a map:
$$\psi:\mathcal{S} \rightarrow \mathbf{C}^{d}.$$
By the Weierstrass preparation 
theorem, the image of $\mathcal{S}$ by $(\psi, Y)$, 
identified in the sequel with 
$\mathcal{S}$, is defined by a unitary polynomial $f \in 
\mathbf{C}\{ X \}[Y]$. The discriminant locus of $\psi$ is defined 
by the discriminant $\Delta_{Y}(f)$ of $f$, which has therefore a d.e.

\begin{defin} \label{qopol}
Let $f \in \mathbf{C}\{ X \}[Y]$ be unitary. If $\Delta_{Y}(f)$ has a 
d.e., we say that $f$ is \textbf{quasi-ordinary}. 
If $\mathcal{A} \simeq \mathbf{C}\{ X \}[Y]/(f)$, with $f$ quasi-ordinary, 
we say that $f$ is a \textbf{qo-defining polynomial} of $\mathcal{S}$ 
and of the algebra $\mathcal{A}$.
\end{defin}

The following theorem, called "of Jung-Abhyankar"  (see \cite{A 55},
\cite{L 88}), generalizes the theorem of Newton-Puiseux for plane
curves. 
 
\begin{theor} \label{JAbh}
If $f \in \mathbf{C}\{ X\}[Y]$ is 
quasi-ordinary, then the set $R(f)$ of roots of $f$ embeds canonically 
in the algebra $\widetilde{\mathbf{C}\{ X\}}$.
\end{theor}
 
In the sequel, we consider $R(f)$ as a subset of 
$\widetilde{\mathbf{C}\{ X\}}$. Moreover, we suppose that $f$ is
\textit{irreducible}.  
Then all the differences of roots of $f$ have d.e., which are totally
ordered for the componentwise order (see \cite{L 65}, 
\cite{L 88}). Denote them by 
$A_{1}<\cdots < A_{G}, \: A_{i}= (A^{1}_{i},...,A^{d}_{i}), \forall 
i \in \{1,...,G\}$. 

\begin{defin} \label{charexp}
We say that the vectors $A_{1},...,A_{G} \in \mathbf{Q}_{+}^{d}$ are 
the \textbf{characteristic exponents} and the monomials
$X^{A_{1}},...,X^{A_{G}}$ are the \textbf{characteristic monomials} of
$f$. 
\end{defin}

After possibly permuting the variables 
$X_{1},...,X_{d}$, we can 
ensure that: 
\begin{equation} \label{lexord}
   (A^{1}_{1},...,A^{1}_{G})\geq_{lex}\cdots \geq_{lex}
        (A^{d}_{1},...,A^{d}_{G}).
\end{equation}
Here $\geq_{lex}$ denotes the lexicographic ordering. In what follows,
we suppose that this condition is always verified.

 \begin{defin} \label{normexp} 
   We say that $f$ is a \textbf{normalized}  qo-defining polynomial of
   $\mathcal{S}$ if (\ref{lexord}) is
   verified and either $A_{1}^{2}\neq 0$ or
   $A_{1}^{1}>1$. 
 \end{defin}

Lipman \cite{L 65} proved that any irreducible quasi-ordinary germ of 
hypersurface has 
normalized qo-defining polynomials (see also \cite{L 83}, \cite{L 88} 
and \cite{GP 00}).

Following Lipman \cite{L 88}, we define inductively the abelian groups
$M=M_{0}:= \mathbf{Z}^{d},\: M_{i}:= M_{i-1}+ \mathbf{Z}A_{i}, \: \forall 
i\in \{1,...,G\}$ and 
the successive indices $N_{i}:= \mbox{card}(M_{i}/M_{i-1}), 
\: \forall 
i\in \{1,...,G\}$. Following Gonz{\'a}lez P{\'e}rez \cite{GP 00}, \cite{GP
  02}, we define 
the vectors $\overline{A}_{1},...,\overline{A}_{G}\in
\mathbf{Q}_{+}^{d}$ as:
\begin{equation} \label{relrec}
    \overline{A}_{1}:= A_{1}, \: 
     \overline{A}_{i}:= N_{i-1} \overline{A}_{i-1} +A_{i}-A_{i-1}, 
          \: \forall  i \in \{2, ...,G\}, \: 
    \overline{A}_{G+1}:= \infty.
\end{equation}

It can be easily seen that $M_{G}$ is also generated by $M_{0}$ and
$\overline{A}_{1},...,\overline{A}_{G}$. Moreover, one has a canonical
way of writing the elements of $M_{G}$ (see \cite{GP 00}, \cite{GP 02}):

 \begin{lemm} \label{uniq}
  Every element of $M_{G}$ can be written in a unique way as a sum $A+
  i_{0} \overline{A}_{1} + \cdots +i_{G-1} \overline{A}_{G}$, where $A
  \in M_{0}$ and $0 \leq i_{k} \leq N_{k+1}-1, \: \forall k \in
  \{0,...,G-1\}$. 
\end{lemm}
 
\textbf{Proof:} As $N_{k}=\mbox{card}(M_{k}/M_{k-1})$, we deduce that $N_{k}
    \overline{A}_{k} \in M_{k-1}, \: \forall k \in \{1,...,G\}$. From
    this we deduce immediately the \textit{existence} of an expression
    verifying the asked property. 

 In order to prove the \textit{uniqueness}, notice first that, 
if $i \in \{1,...,G\}$, one has:
$ N_{i}  = \mbox{min} \{ k \in \mathbf{N}^{*}, \: k A_{i} \in
   M_{i-1}\}
     =\mbox{min} \{ k \in \mathbf{N}^{*}, \: k \overline{A}_{i} \in
   M_{i-1}\}$.  
Then, suppose by contradiction that $\exists (i_{0},...,i_{G-1})\neq
(j_{0},...,j_{G-1})$ and $A,B \in M$ so that $A+ i_{0}
\overline{A}_{1}+\cdots +i_{G-1}\overline{A}_{G}=B+ j_{0}
\overline{A}_{1}+\cdots +j_{G-1}\overline{A}_{G}$. Define $p:=
\mbox{min}\{k \in \{0,...,G-1\},\: i_{l}=j_{l}, \: \forall l >k\}$.
Then $p \geq 0$ and $(i_{p}-j_{p})\overline{A}_{p+1}=(B-A)+
\sum_{k=0}^{p-1}(j_{k}-i_{k})\overline{A}_{k+1}\in M_{p}$.  But
$0<\mid i_{p}-j_{p}\mid \leq N_{p+1}-1$, which contradicts the
previous remark. \hfill $\Box$

\section{The singular locus of the germ} \label{Strucsing}

The characteristic exponents allow to describe precisely the singular locus 
$\mathrm{Sing}(\mathcal{S})$ of $\mathcal{S}$. Before stating this 
description, we introduce some notations taken from \cite{L 88}.

If $I \subset \{1,...,d\}$, let $\mathcal{D}_{I}$ be the linear subspace of 
$\mathbf{C}^{d}$ defined by $\{X_{i}=0, \: \forall i \in I\}$. Its codimension 
is $\mid I \mid$. Denote $\mathcal{Z}_{I}:= \mathcal{S} \cap \psi^{-1}
(\mathcal{D}_{I})$. In 
\cite{L 88}, Lipman shows that the spaces $\mathcal{Z}_{I}$ are
\textit{irreducible}. For 
simplicity, we denote $\mathcal{Z}_{i}:= \mathcal{Z}_{\{i\}}, \: 
\forall i \in \{1,...,d\}$.

\begin{defin} \label{equidim}
The minimal number $c \in \{ 1,...,d\}$ with the property 
that $A_{i}^{k}=0, \: \forall i \in 
\{1,...,G\}, \: \forall k \in \{c+1,...,d\}$ is called 
\textbf{the equisingular dimension} of the quasi-ordinary projection $\psi$.
\end{defin}

 Recalling that the
characteristic exponents verify condition (\ref{lexord}), we see that 
$c$ represents
the number of variables appearing with non-zero exponents among the
monomials $X^{A_{1}},...,X^{A_{G}}$. The name is motivated by the fact
that $\psi$ is then an equisingular deformation of a
$c$-dimensional quasi-ordinary germ, but not of a smaller-dimensional
germ (see Ban \cite{B 93}).

The following theorem is a reformulation of theorem 7.3 in Lipman \cite{L 88}:

\begin{theor} \label{strucsing} 
The irreducible components of $\mathrm{Sing}(\mathcal{S})$ are of the form 
$\mathcal{Z}_{I}$, with $I \subset \{1,...,c\}$ and $\mid I \mid \in \{1,2\}$. 
Moreover:

1) If $i \in \{1,...,c\}$, then $\mathcal{Z}_{i}$ is a component of 
$\mathrm{Sing}(\mathcal{S})$ if and only if one has not simultaneously 
$A_{k}^{i}=0, \: \forall k \in \{1,...,G-1\}$ and $A_{G}^{i}=
\frac{1}{N_{G}}$.

2) If $\{i,j\}\subset \{1,...,c\}$ with $i \neq j$, then 
$\mathcal{Z}_{\{i,j\}}$ is 
a component of $\mathrm{Sing}(\mathcal{S})$ if and only if neither 
$\mathcal{Z}_{i}$ nor $\mathcal{Z}_{j}$ are components of 
$\mathrm{Sing}(\mathcal{S})$.

3) If $\{i,j\}\subset \{1,...,c\}$ with $i \neq j$ and
   $\mathcal{Z}_{\{i,j\}}$ is a component of
   $\mathrm{Sing}(\mathcal{S})$, then the germ of $\mathcal{S}$ at any
   point $P$ of $\mathcal{Z}_{\{i,j\}}-\cup_{k \notin
   \{i,j\}}\mathcal{Z}_{k}$ is isomorphic to the subgerm of
   $(\mathbf{C}^{d+1},0)$ defined by the equation $Y^{N_{G}}=T_{1}T_{2}$.
\end{theor}

\medskip

Let $s \in \{0,...,c\}$ be such that $\mathcal{Z}_{i} \subset
\mathrm{Sing}(\mathcal{S})$ for $i \in \{1,...,s\}$ and $\mathcal{Z}_{i}
\nsubseteq 
\mathrm{Sing}(\mathcal{S})$ either. Our ordering convention
(\ref{lexord}) and the previous theorem imply that there is always
such an $s$, which is equal to \textit{the number of components of
$\mathrm{Sing}(\mathcal{S})$ having codimension 1 in $\mathcal{S}$}. Then 
$\mathcal{Z}_{\{i,j\}}$ is a component of $\mathrm{Sing}(\mathcal{S})$  if and
only if $i \neq j$ and $i,j \in \{s+1,...,c\}$. We get:

\begin{equation} \label{compsing}
\mathrm{Sing}(\mathcal{S})= \bigcup_{1 \leq i \leq s}\mathcal{Z}_{i}
\cup \bigcup_{\stackrel{ j,l \in [[s+1, c]]}{j \neq l}}
\mathcal{Z}_{\{j,l\}}
\end{equation}

\medskip

\noindent \textbf{Examples} (where $G=1$ and $c=3$):

\textbf{1)} If $\mathcal{S}=\{(X,Y),\: Y^{N}-X_{1}X_{2}X_{3}=0\}$, 
with $N >1$, then $s=0$ and $\mathrm{Sing}(\mathcal{S})=
\mathcal{Z}_{\{1,2\}}\cup 
 \mathcal{Z}_{\{2,3\}}\cup \mathcal{Z}_{\{3,1\}}$.

\textbf{2)} If $\mathcal{S}=\{(X,Y),\: Y^{N}-X_{1}^{B^{1}}X_{2}X_{3}=0\}$, 
with $N >1, \: B^{1}>1$, then $s=1$ and  $\mathrm{Sing}(\mathcal{S})=
 \mathcal{Z}_{1}\cup \mathcal{Z}_{\{2,3\}}$.

\textbf{3)} If $\mathcal{S}=\{(X,Y),\: Y^{N}-X_{1}^{B^{1}}X_{2}^{B^{2}}
X_{3}=0\}$, 
with $N >1, \: B ^{1}>1, \: B ^{2}>1$, then $s=2$ and  
$\mathrm{Sing}(\mathcal{S})=\mathcal{Z}_{1}\cup \mathcal{Z}_{2}$.

\textbf{4)} If $\mathcal{S}=\{(X,Y),\: Y^{N}-X_{1}^{B ^{1}}X_{2}^{B ^{2}}
X_{3}^{B ^{3}}=0\}$, 
with $N >1, \: B ^{1}>1, \: B ^{2}>1, \: B ^{3}>1, \: \mbox{gcd}(N,B
^{1}, B ^{2},B ^{3})=1$, then $s=3$ and  
$\mathrm{Sing}(\mathcal{S})=\mathcal{Z}_{1}\cup \mathcal{Z}_{2}\cup
\mathcal{Z}_{3}$.

\section{A reminder of toric geometry} \label{remtor}

For the constructions of sections \ref{norgerm}, \ref{orbmap},
\ref{canoblo}, we need some elementary 
results of toric geometry. For 
details, one can consult Fulton's \cite{F 93} or Oda's \cite{O 88} book.

A \textit{lattice} $\mathcal{W}$ is a finitely generated free abelian
group. An element $v \neq 0$ of $\mathcal{W}$ is called
\textit{primitive} if it cannot be written $v=a v'$ with $a \in
\mathbf{N}^{*}-\{1\}, \: v' \in \mathcal{W}$.  
The dual lattice $\mathcal{M}$ of $\mathcal{W}$ is by definition
$\mbox{Hom}(\mathcal{W}, \mathbf{Z})$. If $\phi:\overline{\mathcal{W}}
\rightarrow \mathcal{W}$ is a 
morphism of lattices and $K$ denotes a field, denote
by $\mathcal{W}_{K}$ the $K$-vector space generated by
$\mathcal{W}$ and by $\phi_{K}$ the corresponding morphism of vector
spaces. 

The elements of $\mathcal{M}$ should be thought about here 
as exponents of \textit{monomials} and those of $\mathcal{W}$ as
\textit{weights} of those monomials.  

Let $\sigma$ be a \textit{strictly convex rational polyhedral}
(abbreviated "s.c.r.p.") \textit{cone} in $\mathcal{W}_{\mathbf{R}}$ and
$\check{\sigma} := \{ u \in \mathcal{M}_{\mathbf{R}}, \: (u,v)\geq 0,
\: \forall v \in \sigma\}$ its dual cone inside
$\mathcal{M}_{\mathbf{R}}$. This dual cone is again s.c.r.p., by Gordan's
lemma. The \textit{affine toric variety} $\mathcal{Z}(\mathcal{W},
\sigma)$ \textit{of weight lattice} $\mathcal{W}$ \textit{and cone} 
$\sigma$ is by definition: 
$$\mathcal{Z}(\mathcal{W},\sigma):=
\mbox{Spec}\mathbf{C}[\check{\sigma} \cap \mathcal{M}]$$
where $\mathbf{C}[\check{\sigma} \cap \mathcal{M}]$  denotes the
$\mathbf{C}$-algebra of the additive semigroup $\check{\sigma} \cap
\mathcal{M}$. Denote by $\chi ^{m}$ the monomial of
$\mathbf{C}[\check{\sigma} \cap \mathcal{M}]$ which has the exponent $m$.

If $v_{1},...,v_{k}\in \mathcal{W}$, we denote by
$\mathbf{R}_{+}(v_{1},...,v_{k})$ the cone generated by them. A
s.c.r.p. cone is 
called \textit{regular} if it is generated by a subset of a basis of 
the lattice $\mathcal{W}$. The variety
$\mathcal{Z}(\mathcal{W},\sigma)$ is smooth if and only if the cone
$\sigma$ is regular.

If $\phi: \overline{\mathcal{W}} \rightarrow \mathcal{W}$ is a
morphism and $\overline{\sigma}\subset
\overline{\mathcal{W}}_{\mathbf{R}}, \: \sigma \subset
\mathcal{W}_{\mathbf{R}}$ are s.c.r.p. cones such that
$\phi_{\mathbf{R}}(\overline{\sigma})\subset \sigma$, 
then there is a canonical
induced \textit{toric morphism} $\phi_{*}: 
\mathcal{Z}( \overline{\mathcal{W}}, \overline{\sigma})\rightarrow 
\mathcal{Z}(\mathcal{W},\sigma)$. 

If $\phi= \mbox{id}_{\mathcal{W}}$,
then $\phi_{*}$ is an embedding. 
As a particular case of this, by taking $\phi=
\mbox{id}_{\mathcal{W}}$  and $\overline{\sigma}=\{0\}$, we obtain a
canonical embedding of the \textit{complex torus} $T_{\mathcal{W}}:=
  \mathcal{Z}( \mathcal{W}, 
\{0\})\simeq (\mathbf{C}^{*})^{\mbox{dim}\mathcal{W}_{\mathbf{R}}}$ in
any affine toric variety 
$\mathcal{Z}(\mathcal{W}, \sigma)$. Moreover,
there is a canonical action of the torus $T_{\mathcal{W}}$ on $\mathcal{Z}(
\mathcal{W}, \sigma)$, such that $\mathcal{Z}(
\mathcal{W},\{0\})$ is the only open orbit. With respect to these
actions, the preceding morphisms are \textit{equivariant}.

If $\phi$ is an inclusion of finite
index and $\phi_{\mathbf{R}}(\overline{\sigma})= \sigma$, 
then $\phi_{*}$ is a finite map. More precisely, we have the following
proposition (corollary 1.16 of \cite{O 88}):

\begin{prop} \label{orb}
If $\phi: \overline{\mathcal{W}} \rightarrow \mathcal{W}$ presents
$\overline{\mathcal{W}}$ as a submodule of finite index of
$\mathcal{W}$, then $\phi_{*}:\mathcal{Z}( \overline{\mathcal{W}},
\overline{\sigma})\rightarrow \mathcal{Z}( \mathcal{W}, \sigma)$ coincides with
the projection for the quotient of $\mathcal{Z}( \overline{\mathcal{W}},
\overline{\sigma})$ with respect to the natural action of the finite group
$\mathcal{W}/ \phi(\overline{\mathcal{W}})$.
\end{prop}

The natural action alluded to comes from the following action of
$\mathcal{W}$ on the monomials of
$\mathbf{C}[\check{\overline{\sigma}}]$:
\begin{equation} \label{actmon}
 v.\chi ^{\overline{u}}=e ^{2i\pi (\overline{u},v)}\chi
 ^{\overline{u}}.
\end{equation}

Let us express the morphism $\phi_{*}$ using coordinates in the case
in which $\mbox{dim } \mathcal{W}_{\mathbf{R}}= \mbox{dim }
\overline{\mathcal{W}}_{\mathbf{R}}=d$ and $\sigma, \overline{\sigma}$
are regular cones of the maximal dimension. Let $v_{1},...,v_{d}\in
\mathcal{W}$ and
$\overline{v}_{1},...,\overline{v}_{d}\in \overline{\mathcal{W}}$ be
the unique primitive elements situated on the edges of $\sigma$,
respectively $\overline{\sigma}$. Write $\phi(\overline{v}_{j}):=
\sum_{i=1}^{d} \alpha_{i}^{j} v_{i}, \: \forall \: j\in
\{1,...,d\}$. The hypothesis $\phi(\overline{\sigma})\subset \sigma$
implies that $\alpha_{i}^{j} \in \mathbf{N}, \: \forall (i,j) \in
\{1,...,d\}^{2}$. Let $u_{1},...,u_{d} \in \mathcal{M}$ and $
\overline{u}_{1},...,\overline{u}_{d}\in \overline{\mathcal{M}}$ be 
the dual bases of   $v_{1},...,v_{d}$, respectively
$\overline{v}_{1},...,\overline{v}_{d}$. The adjoint morphism
$\check{\phi}:\mathcal{M}\rightarrow \overline{\mathcal{M}}$ verifies
then $\check{\phi}(u_{i})= 
\sum^{d}_{j=1}\alpha_{i}^{j} \overline{u}_{j}, \: \forall  i \in
\{1,...,d\}$. The monomials $U_{i}:= \chi ^{u_{i}}, \: 1 \leq i \leq
d$ and $\overline{U}_{i}:= \chi ^{\overline{u}_{i}}, \: 1 \leq i \leq
d$ are free generators of the group algebras
$\mathbf{C}[\check{\sigma}\cap \mathcal{M}]$, respectively
$\mathbf{C}[\check{\overline{\sigma}}\cap \mathcal{M}]$. Then:

\begin{lemm} \label{monex}
The morphism
$\phi ^{*}: \mathbf{C}[\check{\sigma}\cap \mathcal{M}] \rightarrow
\mathbf{C}[\check{\overline{\sigma}}\cap \mathcal{M}] $
can be expressed as: 
$$\begin{cases}
   U_{1}=  \overline{U}_{1}^{\alpha_{1}^{1}}\cdots
  \overline{U}_{d}^{\alpha_{1}^{d}}\\  \vdots \\ 
   U_{d}=  \overline{U}_{1}^{\alpha_{d}^{1}}\cdots
  \overline{U}_{d}^{\alpha_{d}^{d}} 
\end{cases}$$
\end{lemm}

This shows that with respect to the coordinates $U, \overline{U}$ of
the two algebras, the morphism $\phi_{*}$ is \textit{monomial}. Let us look at
its effect on the Newton polyhedron of a fractional series in the
coordinates $U$. If $\eta \in \widetilde{\mathbf{C}\{U\}}, \: \eta =
\sum_{m \in \mbox{Supp}_{U}(\eta)}c_{m}U ^{m}$, where
$\mbox{Supp}_{U}(\eta) \subset 
M_{\mathbf{Q}}$, one has $\phi ^{*}(\eta)= \sum_{m \in
  \mbox{Supp}_{U}(\eta)}c_{m}\overline{U} ^{\check{\phi}(m)}$. So,
\linebreak 
$\mbox{Supp}_{\overline{U}}(\phi ^{*}(\eta)) \subset
\check{\phi}(\mbox{Supp}_{U}(\eta))$. If $m,m' \in M$ and $m' \in
\{m\}+\check{\sigma}$, then from the condition $\phi
^{*}(\check{\sigma})\subset \check{\overline{\sigma}}$ we get
immediately that $\check{\phi}(m')\in \{ \check{\phi}(m)\}+
\check{\overline{\sigma}}$. We deduce:

\begin{lemm} \label{transf}
The vertices of $\mathcal{N}_{\overline{U}}(\phi ^{*}(\eta))$ are
images of vertices of $\mathcal{N}_{U}(\eta)$.
\end{lemm}
\medskip

Define a \textit{fan} $\Sigma$ in $\mathcal{W}_{\mathbf{R}}$ to be a
  finite collection of s.c.r.p. cones, such that for any $\sigma \in
  \Sigma$, all the faces are also in $\Sigma$, and the intersection of
  any two elements of $\Sigma$ is also in $\Sigma$. For example, to
  any s.c.r.p. cone $\sigma$ is associated canonically a fan, the set
  of the faces of $\sigma$. The
  \textit{support} $\mid \Sigma \mid$ of the fan $\Sigma$ is by definition the
  union of the cones composing it. 

Using the
  affine toric morphisms defined before,  
  the affine toric varieties $\mathcal{Z}(
 \mathcal{W}, \sigma), \: \mbox{for } \sigma \in \Sigma$, can be
  equivariantly glued in a new
  variety $\mathcal{Z}(\mathcal{W}, \Sigma)$, called \textit{the toric
  variety of weight lattice $\mathcal{W}$ and fan $\Sigma$}. It
  is always normal. It is smooth if and only if the fan $\Sigma$ is
  \textit{regular}, i.e. its constituting cones are all regular. 

The \textit{orbits} of the action of $T_{\mathcal{W}}$ on
$\mathcal{Z}(\mathcal{W},\Sigma)$ are in 1-to-1 correspondence with
the cones of $\Sigma$. Denote by $O_{\sigma}$ the orbit corresponding
to $\sigma \in \Sigma$ and by $V_{\sigma}$ its closure. One has the
equality $\mbox{dim}\:V_{\sigma}+ \mbox{dim}\:\sigma =
\mbox{dim}\: \mathcal{W}_{\mathbf{R}}$. 

The notion of toric morphism can be extended to this more general
setting, starting from a morphism $\phi: \overline{\mathcal{W}}
\rightarrow \mathcal{W}$ verifying $\forall \: \overline{\sigma} \in
\overline{\Sigma}, \: \exists \sigma \in \Sigma$ such that 
$\phi_{\mathbf{R}}(\overline{\sigma})\subset \sigma$. Then one obtains
an associated (equivariant) toric morphism $\phi_{*}: 
\mathcal{Z}( \overline{\mathcal{W}}, \overline{\Sigma})\rightarrow 
\mathcal{Z}(\mathcal{W},\Sigma)$. A particular case of this
construction is obtained if one considers a \textit{subdivision} of the
fan $\Sigma$, i.e. a second fan $\overline{\Sigma}$ in
$\mathcal{W}_{\mathbf{R}}$ such that $\mid \overline{\Sigma} \mid =
\mid \Sigma \mid $ and $\forall \: \overline{\sigma} \in
\overline{\Sigma}, \: \exists \sigma \in \Sigma$ with 
$\overline{\sigma}\subset \sigma$. The associated toric morphism is 
then proper and birational.

Let us specialize this even more. Consider a regular fan $\Sigma
\subset \mathcal{W}_{\mathbf{R}}$ and let $\sigma_{0}$ be one of its cones, of
dimension $e \geq 2$. Let $V_{1},...,V_{e}$ be the primitive elements
of $\mathcal{W}$ situated on the edges of $\sigma_{0}$ and
$V_{0}:= V_{1}+\cdots + V_{e}$. Denote $\sigma_{0}^{j}:=
\mathbf{R}_{+}(V_{0},..., V_{j-1},V_{j+1},...,V_{e}), \: \forall j \in
\{1,...,e\}$. Each $\sigma \in \Sigma$ with
$\sigma_{0} \subset \sigma$ can be written uniquely as $\sigma=
\sigma_{0}+ \tau$, with  $\tau \in \Sigma$  and $\sigma_{0}\cap \tau
=\{0\}$.  Denote then $\sigma ^{j}:= \sigma_{0}^{j}+\tau, \: \forall j
\in \{1,...,e\}$. 

\begin{defin} \label{star}
 The \textbf{star-subdivision} $\Sigma ^{*}(\sigma_{0})$ of $\Sigma$
 with respect to $\sigma_{0}$ is the fan:
 $ (\Sigma -\{ \sigma \in \Sigma, \: \sigma_{0} \subset \sigma\})
 \cup \{ \mbox{the faces of } \sigma ^{j}, \:  \sigma \in \Sigma,
 \sigma_{0} \subset \sigma, \: 1 \leq j \leq e\}$. In the particular
 case in which $\Sigma$ is the fan 
 composed of the faces of $\sigma_{0}$, we write $\sigma_{0}^{*}$
 instead of $\Sigma ^{*}(\sigma_{0})$.
\end{defin}

The importance of this particular subdivision comes from the
following proposition (see \cite{O 88}), which shows that the
associated toric morphism is intrinsic from an analytical viewpoint:
\begin{prop} \label{torblow}
 The equivariant morphism obtained by passing from $\Sigma$ to the
 star-subdivision $\Sigma ^{*}(\sigma_{0})$ is isomorphic to the
 blow-up of $\mathcal{Z}(\mathcal{W}, \Sigma)$ along
 $V_{\sigma_{0}}$. 
\end{prop}

\section{A toric normalization of the germ} \label{norgerm}

In this section we show how one can obtain using toric geometry a
canonical normalization  
morphism of the germ $(\mathcal{S},0)$ in the chosen coordinate 
system $X$. These are results obtained in \cite{PP 01}.

We denote $W=\mathbf{Z}^{d}$ and let $\sigma_{0}$ be the cone
generated by the canonical basis of $\mathbf{Z}^{d}$. It is a regular
cone. We identify the space $\mathbf{C}^{d}$ of coordinates $X$ with
$\mathcal{Z}(W,\sigma_{0})$. Denote by $u_{1},...,u_{d}$ the
primitive elements of $M=\mbox{Hom}(W, \mathbf{Z})$ situated on the
edges of $\check{\sigma}_{0}$, such that $X_{i}=\chi ^{u_{i}}, \:
\forall i \in \{1,...,d\}$.

Coming back to the quasi-ordinary morphism $\psi$, denote 
$\mathcal{D}:=\mathcal{D}_{1}\cup\cdots \cup \mathcal{D}_{d}$ and 
$\mathcal{Z}:=\mathcal{Z}_{1}\cup\cdots \cup \mathcal{Z}_{d}$. 
Then $\mathcal{S}-\mathcal{Z} \stackrel{\psi}{\rightarrow} 
\mathbf{C}^{d}-\mathcal{D}$ 
is an unramified covering. As the fundamental group 
$\pi_{1}(\mathbf{C}^{d}-\mathcal{D})$ is abelian, this covering is
galoisian, 
its Galois group being $W/W(\psi)$. Here $W(\psi)$ denotes the following 
subsemigroup of $W$:
$$W(\psi):=\psi_{*}\pi_{1}(\mathcal{S}-\mathcal{Z}) \hookrightarrow 
   \pi_{1}(\mathbf{C}^{d}-\mathcal{D})=W.$$

Apart from $\mathbf{C}^{d}=\mathcal{Z}(W, \sigma_{0})$, consider also
the affine $d$-dimensional toric variety $\mathcal{Z}(W(\psi),
\sigma_{0})$ 
and the canonical finite toric morphism:
$$\gamma_{W:W(\psi)}: (\mathcal{Z}(W(\psi), \sigma_{0}),0) 
  \rightarrow (\mathcal{Z}(W, \sigma_{0}),0),$$
where in both cases we denote by $0$ the point which is the unique closed 
orbit of the corresponding toric variety.

By proposition \ref{orb}: 
$\mathcal{Z}(W, \sigma_{0})\simeq \mathcal{Z}(W(\psi), \sigma_{0})/
(W/W(\psi))$. 
The restriction 
$$\gamma_{W:W(\psi)}: \mathcal{Z}(W(\psi), \sigma_{0})- 
        \gamma_{W:W(\psi)}^{-1}(\mathcal{D})
                  \longrightarrow \mathcal{Z}(W, \sigma_{0})-\mathcal{D}$$
is an unramified covering with automorphism group precisely $W/W(\psi)$, which 
shows that one can complete a commutative diagram:
$$\xymatrix{
    \mathcal{Z}(W(\psi), \sigma_{0})-\gamma_{W:W(\psi)}^{-1}(\mathcal{D})
                  \ar[dr]_{\gamma_{W:W(\psi)}} 
            \ar[rr]^-{\nu} & 
           & \mathcal{S}- \mathcal{D} 
                  \ar[dl]^{\psi} \\
    & \mathbf{C}^{d} & }$$
The morphism $\nu$ can be extended to a continuous morphism from 
$(\mathcal{Z}(W(\psi), \sigma_{0}),0)$ to $(\mathcal{S},0)$. As 
the variety $(\mathcal{Z}(W(\psi), \sigma_{0}),0)$ is normal, by the Riemann 
extension theorem, the morphism $\nu$ is in fact everywhere analytic, which 
shows that it is a \textit{normalization} of
$(\mathcal{S},0)$. Noticing that we had not used the fact that
$\mathcal{S}$ is a hypersurface, we get:

\begin{theor} \label{normtor}
For any morphism
$\psi:\mathcal{S}\rightarrow \mathbf{C}^{d}$, unramified over
$(\mathbf{C}^{*})^{d}$ and such that $\mathcal{S}$ is irreducible, one has
the following commutative diagram, in which $\nu$  
is a normalization morphism:
$$\xymatrix{
    (\mathcal{Z}(W(\psi), \sigma_{0}), 0)
         \ar[dr]_{\gamma_{W:W(\psi)}} 
             \ar[rr]^-{\nu} & 
       &  (\mathcal{S}, 0)
                 \ar[dl]^{\psi} \\
    & (\mathbf{C}^{d},0) & }$$
\end{theor}

In the special case in which $\mathcal{S}$ is a hypersurface germ, we
can express the lattice $W(\psi)$ using the characteristic 
exponents of $\psi$. In order to do this, let us introduce the dual 
lattices $W_{k}$ of the lattices $M_{k}$ defined in section \ref{basdef}.
One has the inclusions:
$   M=M_{0} \subsetneq M_{1} \subsetneq \cdots \subsetneq M_{G}, \:
    W=W_{0} \supsetneq W_{1} \supsetneq \cdots \supsetneq W_{G}.$

Introduce also the sequence of fields extensions ( see \cite{L 88} and  
 \cite{GP 00}):
$\mbox{Frac}(\mathbf{C}\{X\})=L=L_{0} \subsetneq L_{1} \subsetneq \cdots 
\subsetneq L_{G}$,  
where:
$L_{k}:= L(X^{A_{1}},...,X^{A_{G}}), \forall k \in \{0,...,G\}$. 
One has the following lemma, proved in \cite{GP 00} (see also \cite{PP 01}):

\begin{lemm} \label{identext}

1) For every $i,j \in \{0,...,G\}, \: i < j$, the fields extension  
$L_{j}:L_{i}$ 
is galoisian and  
$\:\mathrm{Gal}(L_{j}:L_{i})\simeq  W_{j}/ W_{i}.$

2) If $\xi \in R(f)$ then 
$\:\mathrm{Frac}(\mathcal{A})=L(\xi)=L_{G}.$

3) If $N=d_{Y}(f)$, then $\:N=N_{1}\cdots N_{G}$. 
\end{lemm}

The action of the group $W/ W(\psi)=\mbox{Gal}(L(\xi): L)$ on the field 
$L(\xi)$ can be  
canonically lifted as an action of the group $W$
by multiplication with roots of the unity on the monomials of $L(\xi)$
(compare with the relation (\ref{actmon})):
$$v.X^{u}:= e^{2i\pi \langle v, u \rangle}X^{u}, \: \forall v \in W.$$
Here $u$ varies through the set of exponents of the monomials in $L(\xi)$, 
i.e. through the lattice $M_{G}$. So 
$    W(\psi)  = \{v\in W, \: \langle v, u \rangle \in \mathbf{Z}, 
                     \forall u \in M_{G}\}
                = W \cap \mbox{Hom}(M_{G}, \mathbf{Z})= W_{G}$. 
We got like this:

\begin{prop} \label{exprdiff}
Let $f \in \mathbf{C}\{X\}[Y]$ an irreducible quasi-ordinary polynomial and  
$\psi$ be the associated quasi-ordinary projection. Then  $W(\psi)= W_{G}.$
\end{prop}

Using this identification, theorem \ref{normtor} becomes:

\begin{corollary} \label{normtor1}
If $f$ is an irreducible quasi-ordinary polynomial defining the germ
$\mathcal{S}$, then 
one has the following commutative diagram, in which $\nu$ is a normalization  
morphism:
$$\xymatrix{
    (\mathcal{Z}(W_{G}, \sigma_{0}), 0)
         \ar[dr]_{\gamma_{W:W_{G}}} 
             \ar[rr]^-{\nu} & 
       &  (\mathcal{S}, 0)
                 \ar[dl]^{\psi} \\
    & (\mathbf{C}^{d},0) & }$$
\end{corollary}
\medskip

An algebraic proof of this result was obtained by Gonz{\'a}lez P{\'e}rez
in \cite{GP 00}. In the case of surfaces, the normalization of an
irreducible quasi-ordinary germ has a Jung-Hirzebruch singularity (see
Barth-Peters-Van de Ven \cite{BPV 84}).

\section{The canonical orbifold map} \label{orbmap}

In this section we show that there is a canonical finite morphism 
$\mu$ whose target is the normalization $\overline{\mathcal{S}}$ of 
$\mathcal{S}$, its source being a smooth germ $\mathcal{R}$. It is 
a particular case of the \textit{orbifold maps} defined by Deligne and
Mostow in \cite{DM 93}, and 
also of a construction described by Prill in \cite{P 67} using
Grauert-Remmert's existence theorem (see Bell-Narasimhan \cite{BN 90}
for a presentation of this last theorem). 

 We denote again by $0$ the base point $\nu ^{-1}(0)$ of
 $\overline{\mathcal{S}}$.  
Consider the extrinsic isomorphism $\overline{\mathcal{S}} 
\simeq \mathcal{Z}(W_{G}, \sigma_{0})$ of the previous section. Let 
$\tilde{W}$ be the sublattice of $W_{G}$ generated by the smallest non-zero  
elements of $W_{G}$ situated on the edges of $\sigma_{0}$. Then 
$\sigma_{0}$ is regular with respect to $\tilde{W}$, and so 
$\mathcal{Z}(\tilde{W}, \sigma_{0})$ is smooth. Consider the toric map:
$$\mu: \mathcal{Z}(\tilde{W}, \sigma_{0}) \rightarrow 
   \mathcal{Z}(W_{G}, \sigma_{0})$$
obtained by changing the lattice. Denote by
$\tilde{u}_{1},...,\tilde{u}_{d}$ the primitive elements of
$\tilde{M}=\mbox{Hom}(\tilde{W}, \mathbf{Z})$ situated on the edges of
$\check{\sigma}_{0}$, such that the image of $u_{i}$ is proportional
with $u_{i}$ in $M_{\mathbf{R}}$.  Introduce $\tilde{U}_{i}:= \chi
^{\tilde{u}_{i}}, \: \forall i \in \{1,...,d\}$.  Then, the
composition $\nu \circ \mu : \mathcal{Z}(\tilde{W},
\sigma_{0})\rightarrow \mathcal{Z}(W, \sigma_{0})$ is given in the
coordinates $(\tilde{U}_{i})_{1\leq i \leq d}, \: (X_{i})_{1\leq i
  \leq d}$ by equations of the form $X_{i}= \tilde{U}_{i}^{m_{i}}, \:
m_{i}\in \mathbf{N}^{*}, \: \forall i \in \{1,...,d\}$. Denote also by 
$\tilde{v}_{1},...,\tilde{v}_{d}$ the dual basis of
$\tilde{u}_{1},...,\tilde{u}_{d}$.

By proposition \ref{orb}, the morphism $\mu$ is 
the quotient map of $\mathcal{Z}(\tilde{W},
\sigma_{0})$ by the natural action of the
finite group $W_{G}/\tilde{W}$. Moreover, it can be easily seen that
$W_{G}/\tilde{W}$ does not contain  
complex reflections. This shows that the locus $F_{\mu}$ of the fixed
points of the elements of $W_{G}/\tilde{W}$ distinct from the identity has
codimension at least 2 in $\mathcal{Z}(\tilde{W}, \sigma_{0})$. Moreover, 
$\mu^{-1}(\mathrm{Sing}(\mathcal{Z}(W_{G}, \sigma_{0})))\subset
F_{\mu}$.  As $\mathcal{Z}(\tilde{W}, \sigma_{0})$ is smooth, the
complement $\mathcal{Z}(\tilde{W}, \sigma_{0})-  
  \mu^{-1}(\mathrm{Sing}(\mathcal{Z}(W_{G}, \sigma_{0})))$ is 
simply connected, and so the restriction of $\mu$ over the smooth part 
of $\mathcal{Z}(W_{G}, \sigma_{0})$ is the universal covering
map. This shows its uniqueness, by the same arguments as in 
\cite{P 67} or \cite{DM 93}. More precisely, we have the following
result, which is a particular case of proposition 14.3 of \cite{DM
  93}, proved there algebraically:

\begin{lemm} 
For $i =1,2$, let $\mu_{i}: (\mathcal{R}_{i}, P_{i}) \rightarrow
(\overline{\mathcal{S}},0)$ be two finite maps unramified in
codimension 1, with 
smooth sources $\mathcal{R}_{i}$. Then there exists an isomorphism
$(\mathcal{R}_{1}, P_{1})\rightarrow(\mathcal{R}_{2}, P_{2})$ making
the following diagram commutative:
$$\xymatrix{
    (\mathcal{R}_{1},P_{1})
         \ar[dr]_{\mu_{1}} 
             \ar[rr] & 
       &   (\mathcal{R}_{2},P_{2})
                 \ar[dl]^{\mu_{2}} \\
    & (\overline{\mathcal{S}},0) & }$$
\end{lemm}

We denote $\mathcal{R}:= \mathcal{Z}(\tilde{W}, \sigma_{0})$. The preceding 
lemma shows that the morphism $\mu: \mathcal{R} 
\rightarrow \overline{\mathcal{S}}$ is independent  of the particular 
isomorphism $\overline{\mathcal{S}} \simeq \mathcal{Z}(W_{G},
\sigma_{0})$ under consideration.
 
\begin{defin} \label{canorb}
We call the mapping $\mu$ described before \textbf{the canonical
  orbifold mapping} associated to $\overline{\mathcal{S}}$.
\end{defin}

The vocabulary is motivated by the fact that $\mu$
is locally the quotient map of the action of a finite group on $\mathcal{R}$. 

Denote $P_{0}:= (\nu \circ \mu)^{-1}(0)$. So we have constructed a
canonical finite map of germs of analytical spaces $\nu \circ \mu : (\mathcal{R},
P_{0}) \rightarrow (\mathcal{S},0)$.

\section{Expansions according to semiroots} \label{exprts}

Semiroots were introduced for arbitrary irreducible
quasi-ordinary polynomials  
by Gonz{\'a}lez P{\'e}rez in \cite{GP 00}. For their use in the study of plane
curves, see \cite{PP 99}. Here we need them as an essential tool in
the proofs of theorem \ref{isosgsg}.

\begin{defin} \label{defsrac}
Let us fix $\xi \in R(f)$. Take any $k \in \{0,...,G\}$. A unitary polynomial 
$f_{k} \in \mathbf{C}\{X\}[Y]$ is called a 
\textbf{$\mathbf{k}-$semiroot} of $f$ if  
$f_{k}$ is of degree $N_{1}\cdots N_{k}$, and $f_{k}$ has a d.e.  
with $v_{X}(f_{k}(\xi)) = \overline{A}_{k+1}$.  
A $(G+1)$-tuple $(f_{0},...,f_{G})$ such that $\forall k \in \{0,...,G\}, \: 
f_{k}$ is a $k$-semiroot of $f$ is called a \textbf{complete system
of semiroots} for $f$. 
\end{defin}

These objects are independent of the choice of $\xi$.

Let $(f_{0},...,f_{G})$ be a complete system of semiroots for $f$ (which 
always exists, for example the 
minimal polynomials of suitable truncations of $\xi$ or the
characteristic approximate roots of $f$, see \cite{GP 00}). 
Generalizing immediately Abhyankar \cite{A 77} (see also \cite{PP
  99}), we have:

\begin{lemm} 
Any  
$h\in \mathbf{C}\{X\}[Y]$ can be uniquely written as a finite sum 
$h = \sum c_{i_{0}\cdots i_{G}} f_{0}^{i_{0}}\cdots f_{G}^{i_{G}}$,
with $c_{i_{0}\cdots i_{G}}\in \mathbf{C}\{X\}$, 
the  $(G+1)$-tuples 
$(i_{0},...,i_{G})\in \mathbf{N}^{G+1}$ verifying  
$0\leq i_{k} \leq N_{k+1}-1, \: \forall k \in \{0,...,G-1\}$ and
$i_{G}\leq [\frac{d_{Y}(h)}{N}]$. 
\end{lemm}

\textbf{Proof:} Make the euclidean division of $h$ by $f_{G}$ and of
the successive quotients by $f_{G}$, till one obtains a quotient of
degree $< d_{Y}(f_{G})$. Then one gets the $f_{G}$-\textit{adic
  expansion} of $h$ by $f_{G}$, which has the form $h=\sum c_{i_{G}}
f_{G}^{i_{G}}$, where $i_{G} \leq [\frac{d(h)}{d(f_{G})}]$. Then
iterate this, making at each step the $f_{k-1}$-adic expansions of the
coefficients $c_{i_{k}\cdots i_{G}}$. 

The unicity comes from the
remark that the $Y$-degrees of the terms \linebreak$c_{i_{0}\cdots i_{G}}
f_{0}^{i_{0}}\cdots f_{G}^{i_{G}}$ are pairwise distinct. To see it,
suppose by contradiction that $\exists (i_{0},...,i_{G})\neq
(j_{0},...,j_{G})$ and $d_{Y}(c_{i_{0}\dots
  i_{G}}f_{0}^{i_{0}}\cdots f_{G}^{i_{G}})= d_{Y}(c_{j_{0}\dots
  j_{G}}f_{0}^{j_{0}}\cdots f_{G}^{j_{G}})$. Then $\exists p \in
\{0,...,G\}$ so that $i_{k}=j_{k}, \forall k >p$ and $i_{p}\neq
j_{p}$. Suppose for example that $i_{p} >j_{p}$. Then $(i_{p}
-j_{p})N_{1}\cdots N_{p}= \sum_{k=0}^{p-1}(j_{k}-i_{k})N_{1}\cdots
N_{k}\leq \sum_{k=0}^{p-1}(N_{k+1}-1)N_{1}\cdots N_{k}= N_{1}\cdots
N_{p}-1$, and so $i_{p}-j_{p}<1$, which is a contradiction. \hfill
$\Box$

\begin{defin} \label{dvptadq}
The preceding equality is called 
\index{d{\'e}veloppement -adique}\textbf{the 
$(f_{0},...,f_{G})$-adic expansion} of $h$. 
The finite set 
$\{ (i_{0},...,i_{G}), \: c_{i_{0}\cdots i_{G}} \neq 0 \}$ is called  
\index{support}\textbf{the $(f_{0},...,f_{G})$-adic support of $h$}, 
denoted $\mathrm{Supp}_{(f_{0},...,f_{G})}(h)$.
\end{defin}

 The following lemma, which
generalizes the properties of Abhyankar's  
expansions in terms of semiroots in the plane branch case (see \cite{A 77} 
and \cite{PP 99}), is a simple
consequence of lemma \ref{uniq}:

\begin{lemm} \label{corlemm}
If $h = \sum c_{i_{0}\cdots i_{G}} f_{0}^{i_{0}}\cdots f_{G}^{i_{G}}$ 
is the $(f_{0},...,f_{G})$-adic expansion of $h \in 
\mathbf{C}\{X\}[Y]$, then for every $\xi \in R(f)$, the sets of vertices of 
the Newton polyhedra $\mathcal{N}_{X}
(c_{i_{0}\cdots i_{G}} (f_{0}(\xi))^{i_{0}}\cdots (f_{G}(\xi))^{i_{G}})$ 
are pairwise disjoint, when $(i_{0},...,i_{G})$ varies through the 
$(f_{0},...,f_{G})$-adic support of $h$.
\end{lemm}

In the sequel, the previous lemma will be important combined with the
following one:

\begin{lemm} \label{disjpol}
If $h_{1},...,h_{p}\in \widetilde{\mathbf{C}\{X\}}$ and the sets of
vertices of the Newton polyhedra
$\mathcal{N}_{X}(h_{1}),...,\mathcal{N}_{X}(h_{p})$ are pairwise
disjoint, then $\mathcal{N}_{X}(h_{1}+\cdots + h_{p})$ is the convex
hull of the union $\mathcal{N}_{X}(h_{1})\cup \cdots \cup
\mathcal{N}_{X}(h_{p})$. In particular, each vertex of
$\mathcal{N}_{X}(h_{1}+\cdots + h_{p})$ is a vertex of one of the
polyhedra $\mathcal{N}_{X}(h_{1}),...,\mathcal{N}_{X}(h_{p})$. 
\end{lemm}

\textbf{Proof:} Denote by $\mbox{Conv}(E)$ the convex hull of a set
$E$. If  $h:= h_{1}+\cdots +h_{p}$, one has always the inclusion:
$\mathcal{N}_{X}(h)\subset \mbox{Conv}(\mathcal{N}_{X}(h_{1})\cup \cdots \cup
\mathcal{N}_{X}(h_{p}))$,  
as each monomial of $h$ is a monomial of one of the $h_{i}$. 

Conversely, each vertex of $\mbox{Conv}(\mathcal{N}_{X}(h_{1})\cup \cdots \cup
\mathcal{N}_{X}(h_{p}))$ is a vertex of one of the
$\mathcal{N}_{X}(h_{i})$. The hypothesis of the lemma implies that it
is necessarily also a vertex of $\mathcal{N}_{X}(h)$, which proves the
converse inclusion. \hfill $\Box$

\section{Various definitions of semigroups} \label{extrsg}

A difficulty for extending the second definition of the semigroup of a
plane irreducible curve is 
that in dimension $>1$, fractional series may have no dominating exponent. 

In \cite{GP 00}, \cite{GP 02}, 
Gonz{\'a}lez P{\'e}rez considers a more general notion, by using instead of
the dominating exponents, the Newton polyhedra 
$\mathcal{N}_{X}(h(\xi))$ for varying $h \in \mathbf{C}\{X\}[Y]-(f)$  
and the set of their vertices:
$$ \Gamma_{\mathcal{N}}(f):= \{ A \in  \mathbf{Q}_{+}^{d}, \: 
    A \mbox{ is a vertex of } \mathcal{N}_{X}(h(\xi)), \: 
    h \in \mathbf{C}\{X\}[Y]-(f) \}.$$
Here $\xi$ denotes again an arbitrary root of $f$. 
%As $\mathcal{N}_{X}(h_{1}(\xi)h_{2}(\xi))$ is the Minkowski sum
%$\mathcal{N}_{X}(h_{1}(\xi))+ \mathcal{N}_{X}(h_{2}(\xi))$ and each vertex
%of the Minkowski sum of two polyhedra  is the sum of 
%the vertices of the initial polyhedra, we see that
%$\Gamma_{\mathcal{N}}(f)$ is a sub-semigroup of $(\mathbf{Q}_{+}^{d},
%+)$. 

One sees immediately the independence of the set 
$\Gamma_{\mathcal{N}}(f)$ from the choice of
$\xi$. In \cite{GP 00}, \cite{GP 02}, Gonz{\'a}lez P{\'e}rez proves: 

\begin{prop} \label{gensgr}
  One has the equality of sets:
  $$\Gamma_{\mathcal{N}}(f)= \mathbf{N}^{d}+ \mathbf{N} \overline{A}_{1}
  + \cdots + \mathbf{N} \overline{A}_{G}.$$
\end{prop}

\textbf{Proof:} This is an immediate consequence of lemmas
\ref{corlemm} and \ref{disjpol}. \hfill $\Box$

\medskip

So, the set $ \Gamma_{\mathcal{N}}(f)$ has the structure of
semigroup for the addition, a fact which was not a priori clear from
the definition. 

Now we modify the previous definition, by considering only those 
functions which have a dominating exponent, as we find it convenient
for stating  theorem \ref{isosgsg}. 
We introduce the following subset of $\mathbf{Q}_{+}^{d}$:
$$ \Gamma_{\mathcal{D}}(f):= \{v_{X}(h(\xi)), \: h \in
 \mathbf{C}\{X\}[Y]-(f), \: h(\xi) 
\mbox{\textrm{ has a d.e.}}\}.$$
This time it is clear that it is a semigroup. 
In fact one obtains the same semigroup as before, again an immediate
consequence of lemmas \ref{corlemm} and \ref{disjpol}:

\begin{prop} \label{samesgr}
One has the equality of semigroups:
  $$\Gamma_{\mathcal{D}}(f)= \mathbf{N}^{d}+ \mathbf{N} \overline{A}_{1}
  + \cdots + \mathbf{N} \overline{A}_{G}.$$
\end{prop}

This motivates the following definition:

\begin{defin} \label{semextr}
The semigroup $\Gamma_{\mathcal{N}}(f)=\Gamma_{\mathcal{D}}(f)$ is
called \textbf{the semigroup of} $\mathcal{A}$ \textbf{with respect to} $f$, 
denoted $\Gamma(f)$.
\end{defin}

Let us introduce now some notions needed to give in sections \ref{simple} and
\ref{results} another definition of a semigroup associated to $f$,
generalizing the first one of the introduction.

Let $\mathcal{V}$ be a complex analytical manifold of dimension $n$ and 
$\mathcal{H}$ a hypersurface of $ \mathcal{V}$. Let $P \in \mathcal{H}$ 
be a point such that the germ of $\mathcal{H}$ at $P$ has normal crossings. 
Let $(\mathcal{H}_{1},P),..., (\mathcal{H}_{r},P)$ be the irreducible 
components of $(\mathcal{H},P)$, where $1 \leq r \leq n$. From now on,
their ordering is supposed to be fixed.

\begin{defin} \label{centre}
The $r$-codimensional submanifold $C(\mathcal{H},P):= 
(\mathcal{H}_{1},P) \cap \cdots \cap (\mathcal{H}_{r},P)$ of $\mathcal{V}$ 
is called the \textbf{center} of  $\mathcal{H}$ at $P$.
\end{defin}

Let $(x_{1},...,x_{n})$ be local coordinates of $\mathcal{V}$ at $P$, 
such that $\mathcal{H}_{i}=\mathcal{Z}(x_{i}), \: \forall i \in
\{1,...,r\}$. 
We say in this case that they are \textit{adapted} to $\mathcal{H}$ at
$P$.  Take $h \in \mathcal{O}_{\mathcal{V},P}$. 
Write $h=x_{1}^{m_{1}}\cdots x_{r}^{m_{r}}u$, with $u \in \mathcal{O}_
{\mathcal{V},P}$ not divisible by any of $x_{1},...,x_{r}$. As $m_{i}$
is the multiplicity of 
$\mathcal{H}_{i}$ as an irreducible component of the principal divisor
$(h)$, we have:

\begin{lemm} \label{invexp}
The $r$-tuple $(m_{1},...,m_{r})$ depends only on the ordering of the 
components of the germ $(\mathcal{H},P)$, 
and not on the choice of the adapted coordinates $x_{1},...,x_{n}$.
\end{lemm}

Let us define a subclass of elements of
$\mathcal{O}_{\mathcal{V},P}$, distinguished by their relation with
$\mathcal{H}$:

\begin{defin} \label{domexp}
We say that $h \in \mathcal{O}_{\mathcal{V},P}$ \textbf{has a
  dominating exponent (abbreviated "d.e.") 
with respect to the germ} $(\mathcal{H},P)$ if 
$C(\mathcal{H},P)$ is not included in the closure of
  $\mathcal{Z}(h)-\mathcal{H}$.  
The $r$-tuple $(m_{1},...,m_{r})
\in \mathbf{N}^{r}$, written $v_{\mathcal{H},P}(h)$, is called 
\textbf{the dominating exponent} of $h$ with respect to $(\mathcal{H},P)$.
\end{defin}

Lemma \ref{invexp} shows that the dominating exponent is independent
of the chosen adapted coordinates, once we have chosen an ordering of
the components of $\mathcal{H}$ at $P$. 
In coordinates $(x_{1},...,x_{n})$ adapted to $\mathcal{H}$ at
$P$, we see
that $h$ has a d.e. with respect to $(\mathcal{H}, P)$ 
if and only if the series $h \in \mathbf{C}\{x_{1},...,x_{n}\}$
has a d.e. with respect to $(x_{1},...,x_{r})$, in the sense of
section \ref{NP}, i.e. if and only if $u(x_{1},...,x_{r},0,...,0)\neq0$.

If $\mathcal{B}$ is a subalgebra of $\mathcal{O}_{\mathcal{V},P}$, let us 
introduce the following subsemigroup of the multiplicative semigroup 
$(\mathcal{B},\cdot)$:
$$\mathcal{E}_{\mathcal{H},P}(\mathcal{B}):= 
  \{ h \in \mathcal{B}-\{0\}, \: h \: \mbox{ has a d.e. with 
      respect to } \:(\mathcal{H},P)\: \}.$$

This allows us to  define:

\begin{defin} \label{semhyp}
The \textbf{semigroup of $\mathcal{B}$ with respect to} $(\mathcal{H},P)$ 
is the following subsemigroup of $(\mathbf{N}^{r},+)$:
$$\Gamma_{\mathcal{H},P}(\mathcal{B}):= 
      \{ v_{\mathcal{H},P}(h), \: h\in 
       \mathcal{E}_{\mathcal{H},P}(\mathcal{B})-\{0\} \}.$$
\end{defin}

This semigroup consists simply of the dominating exponents with respect
to $(\mathcal{H},P)$ of the elements of $\mathcal{B}$ which have d.e. 

\medskip

\section{The simplest case of the main theorem} \label{simple}

Let us suppose in this section that $s=d$. Then
$\mathrm{Sing}(\mathcal{S})= \bigcup_{1 \leq i \leq d}\mathcal{Z}_{i}$
(see section \ref{Strucsing}). 

Let $\nu:\overline{\mathcal{S}}\rightarrow \mathcal{S}$ be the
normalization morphism of $\mathcal{S}$ and $\mu:
\mathcal{R}\rightarrow \overline{\mathcal{S}}$ be the orbifold map of
$\overline{\mathcal{S}}$. Denote $\theta:=\nu \circ \mu: \mathcal{R}\rightarrow
\mathcal{S}$ and let $\theta ^{*}$ be the corresponding morphism of
sheaves of local algebras. Define $\overline{\mathcal{H}}:=\theta
^{-1}(\mathrm{Sing}(\mathcal{S}))$. Then, using the toric
presentations of the morphisms $\nu$ and $\mu$ given in the sections
\ref{norgerm} and \ref{orbmap}, we see that $\overline{\mathcal{H}}$
is a divisor with normal crossings whose center at $P:=\theta^{-1}(0)$
is 0-dimensional, reduced to $P$ itself. Denote:
$$\Gamma'_{P}(\mathcal{S}):=\Gamma_{\overline{\mathcal{H}},P}(\theta
^{*}(\mathcal{A})_{P}).$$ 
This sub-semigroup of $(\mathbf{N}^{d},+)$ is obviously an
\textit{analytical invariant} of $\mathcal{S}$. The following theorem,
which is the main one of this article specialized to the case $s=d$,
shows that this semigroup is isomorphic to $\Gamma(f)$. The general
case is stated in theorem \ref{isosgsg}.

\begin{theor} \label{mths}
Let $f$ be a quasi-ordinary defining polynomial of
$\mathcal{S}$. Suppose that $s=d$. If $h$ varies through
$\mathbf{C}\{X\}[Y]$ such that $h(\xi)$ has a d.e., then $\theta
^{*}(h \mid_{\mathcal{S}})$ has a d.e. with respect to
$\overline{\mathcal{H}}$ at $P:=\theta^{-1}(0)$ and one obtains a
well-defined mapping:
$$ \begin{array}{lll}\Phi_{P}: &
              \Gamma(f) & \longrightarrow \Gamma'_{P}(\mathcal{S})\\
              & v_{X}(h(\xi)) & \longrightarrow
              v_{\overline{\mathcal{H}},P}(\theta^{*}
                 (h\mid_{\mathcal{S}})_{P})
            \end{array}$$
which realizes an isomorphism of semigroups.
\end{theor}

\textbf{Proof:} Denote $\overline{\psi}:= \psi \circ \theta:
\mathcal{R} \rightarrow  \mathbf{C}^{d}$. 

As the image of $Y\in \mathbf{C}\{X\}[Y]$ in $\mathcal{A}$ verifies the 
equation $f(X, Y)=0$, one sees that $Y\mid_{\mathcal{S}}$ can 
be thought as an element of $R(f)$. Denoting it by $\xi$, one has the 
equality (\ref{relfund}):
$ \: \overline{\psi}^{*}(h(\xi))=\theta^{*}(h\mid_{\mathcal{S}}), \: 
   \forall h \in \mathbf{C}\{X\}[Y].$

Taking the toric representatives constructed
before of the morphisms $\mu$ and $\nu$,  the 
point $P$ is the orbit of dimension $0$ and one can choose 
canonical toric coordinates adapted to $\overline{\mathcal{H}}$ at $P$. 
With such coordinates, the morphism $\overline{\psi}^{*}$ is 
monomial, and using formula (\ref{relfund}), we see that $\Phi_{P}$ is 
\textit{is injective}.

In order to prove its surjectivity, we must show that if 
$h \in \mathbf{C}\{X\}[Y]$ is such that $\theta ^{*}(h\mid_{\mathcal{S}}) \in 
\mathcal{E}_{\overline{\mathcal{H}},P}(\theta ^{*}(\mathcal{A}))$,
then one can find another element  
$h' \in \mathbf{C}\{X\}[Y]$  such that $h(\xi)$ has a d.e. and
$v_{\overline{\mathcal{H}},P}(\theta^{*} 
(h'\mid_{\mathcal{S}})_{P})= v_{\overline{\mathcal{H}},P}(\theta^{*}
(h\mid_{\mathcal{S}})_{P})$. As $f\mid_{\mathcal{S}}=0$, we can suppose 
that $\mbox{deg}(h)< \mbox{deg}(f)$, after possibly making the euclidian 
division of $h$ by $f$. We consider then a complete system 
$(f_{0},...,f_{G})$ of semiroots of $f$ and the $(f_{0},...,f_{G})$-adic 
expansion of $h$, which by our hypothesis is of the form 
$h= \sum c_{i_{0}\cdots i_{G-1}} f_{0}^{i_{0}}\cdots f_{G-1}^{i_{G-1}}$. 
Using lemma \ref{corlemm}, we see that there \linebreak exists a tuple  
$(i_{0},...,i_{G-1},0)
\in \mbox{Supp}_{(f_{0},...,f_{G})}(h)$ such that 
$v_{\overline{\mathcal{H}},P}(\theta^{*}(h\mid_{\mathcal{S}})_{P})= 
v_{\overline{\mathcal{H}},P}(\theta^{*}
(c_{i_{0}\cdots i_{G-1}} f_{0}^{i_{0}}\cdots f_{G-1}^{i_{G-1}}
\mid_{\mathcal{S}})_{P})=
v_{\overline{\mathcal{H}},P}(\theta^{*}
(X^{m} f_{0}^{i_{0}}\cdots f_{G-1}^{i_{G-1}}
\mid_{\mathcal{S}})_{P})$, where $m$ is one of the vertices of the 
Newton polyhedron 
$\mathcal{N}_{X}(c_{i_{0}\cdots i_{G-1}})$. But the term 
$X^{m} f_{0}^{i_{0}}\cdots f_{G-1}^{i_{G-1}}$ has a d.e, which 
proves that $\Phi_{P}$ \textit{is  surjective}.
\hfill $\Box$

\section{A canonical sequence of blowing-ups} \label{canoblo}

In this section we consider $(\nu \circ
\mu)^{-1}(\mathrm{Sing}(\mathcal{S}))$ as a subspace of $\mathcal{R}$
and we construct from it a canonical map $\eta: \overline{\mathcal{R}} 
\rightarrow \mathcal{R}$ obtained as a composition of blowing-ups of
smooth centers. This construction will be used in order to replace the
morphism $\theta$ of theorem \ref{mths} by $\nu \circ \mu \circ \eta$
(see section \ref{results}).

Define:
 $$c' := \begin{cases}
          c-2, \:  \mathrm{ if }\:  s = c-2, \\
          c, \: \mathrm{ if }\: s \neq c-2
  \end{cases} $$
$$ \mathcal{D}':= \mathcal{D}_{1}\cup \cdots \cup \mathcal{D}_{c'}$$
$$   \mathcal{Z}':= \begin{cases}
         \bigcup_{1 \leq i \leq s}\mathcal{Z}_{i} , \:  
               \mathrm{ if }\:  s \in \{c-2, c-1\}, \\
        \bigcup_{1 \leq i \leq s}\mathcal{Z}_{i}\cup \mathcal{Z}_{[[s+1,c]]}, 
                    \: \mathrm{ if }\: s \notin \{c-2,c-1\} \\
         \end{cases} $$
We say that $c'$ is the \textit{reduced equisingular dimension} of
$\psi$ (see also definition \ref{equidim}). When $s \in \{c-1,c\}$,
the space $\mathcal{Z}'$ is  
precisely the union of the
components of $ \mathrm{Sing}(\mathcal{S})$ which have codimension 1
in $\mathcal{S}$. As $\mathcal{Z}_{[[s+1,c]]}= 
\bigcap_{\stackrel{ j,l \in [[s+1,c]]}{j \neq l}}
                   \mathcal{Z}_{\{j,l\}}$, using formula (\ref{compsing})
we see that $\mathcal{Z}'\hookrightarrow \mathrm{Sing}(\mathcal{S})$
depends only on the
analytical structure of $\mathcal{S}$.

When constructing $\eta$, we consider several cases, according to the
 values of $s$, the number 
 of components of $\mathrm{Sing}(\mathcal{S})$ of codimension 1 in
 $\mathcal{S}$, and of the equisingular dimension $c$:

1) \textbf{Suppose that $s \leq c-3$}.

If $I \subset \{1,...,d\}$, denote 
$\mathcal{U}_{I}:=(\nu \circ \mu)^{-1}(\mathcal{Z}_{I})$. Using
the toric construction of the composition $\nu \circ \mu$ presented in
the last two sections, one sees
that in the canonical toric coordinates $\tilde{U}_{1},...,\tilde{U}_{d}$ of
$\mathcal{R}$, one has $\mathcal{U}_{I}=\bigcap_{i \in
  I}\mathcal{Z}(\tilde{U}_{i})$. Formula (\ref{compsing}) shows that:
\begin{equation} \label{preimsing}
      (\nu \circ \mu)^{-1}(\mathrm{Sing}(\mathcal{S}))=
      \bigcup_{1 \leq i \leq s}\mathcal{U}_{i}
      \cup \bigcup_{\stackrel{ j,l \in [[s+1, c]]}{j \neq l}}
      \mathcal{U}_{\{j,l\}}.
\end{equation}

\textbf{Consider first the case in which $c=d$ and $s=0$}. Then
$\mathrm{Sing}(\mathcal{S})$ has only components of codimension 2 in
$\mathcal{S}$: 
$ \mathrm{Sing}(\mathcal{S})=\bigcap_{\stackrel{I \subset
    \{1,...,d\}}{\mid I \mid =2}}\mathcal{Z}_{I}$. 

Then the axis of the coordinates $\tilde{U}_{1},...,\tilde{U}_{d}$ can
be obtained 
analytically from $ (\nu \circ
\mu)^{-1}(\mathrm{Sing}(\mathcal{S}))$. Indeed, if $\mathcal{L}^{i}$
is the axis of the coordinate $\tilde{U}_{i}$, one has 
$\mathcal{L}^{i}= \bigcap_{\stackrel{\mid I \mid =2}{i \notin
    I}}\mathcal{U}_{I}$.

Let 
$\pi_{1}: \mathcal{R}_{1} \rightarrow \mathcal{R}$ be the blow-up of 
$\mathcal{R}$ at $P_{0}= (\nu \circ \mu)^{-1}(0)$. Denote by 
$\mathcal{H}_{1}:= \pi_{1}^{-1}(P_{0})$ the exceptional divisor of
$\pi_{1}$ and by $P_{1}^{i}$ the point where the strict transform of
$\mathcal{L}^{i}$ meets the exceptional divisor $\mathcal{H}_{1}$. 

 Let $\pi_{2}:\mathcal{R}_{2} \rightarrow \mathcal{R}_{1}$
be the blow-up of $\mathcal{R}_{1}$ at all the points $P_{1}^{i}$, for $i
\in \{1,...,d\}$. Denote by
$\mathcal{H}_{2}^{i}:=\pi_{2}^{-1}(P_{1}^{i})$ the components of the
exceptional divisor of $\pi_{2}$, by $\mathcal{H}_{1,2}$ the strict
transform of $\mathcal{H}_{1}$ by $\pi_{2}$.  Denote by 
$\mathcal{L}^{\{i,j\}}$ the line of the $(c-1)$-dimensional
projective space $\mathcal{H}_{1}$ joining $P_{1}^{i}$ and
$P_{1}^{j}$.
Define $P_{2}^{i,j}$ to be the
point where the strict transform of $\mathcal{L}^{\{i,j\}}$ meets the
divisor $\mathcal{H}_{2}^{i}$.

More generally, suppose that $\mathcal{R}_{k-1}$ is already
constructed, with $d \geq
k \geq 2$. Let $\pi_{k}: \mathcal{R}_{k} \rightarrow
\mathcal{R}_{k-1}$ be the blow-up of $\mathcal{R}_{k-1}$ at all the
points $P_{k-1}^{i_{1},...,i_{k-1}}$, for $i_{1},...,i_{k-1}\in 
\{1,...,d\}$ pairwise distinct. The components of the exceptional
divisor of $\pi_{k}$ are $\mathcal{H}_{k}^{i_{1},...,i_{k-1}}:= \pi_{k}^{-1}
(P_{k-1}^{i_{1},...,i_{k-1}})$. For $l \in \{1,...,k-1\}$ and
$j_{1},...,j_{l-1}\in \{1,...,d\}$ pairwise distinct, the strict
transform by $\pi_{k}$ of $\mathcal{H}_{l,k-1}^{j_{1},...,j_{l-1}}$ is
denoted $\mathcal{H}_{l,k}^{j_{1},...,j_{l-1}}$. Denote by
$\mathcal{L}^{i_{1},...,i_{k-2},\{i_{k-1},i_{k}\}}$ the line 
joining the points $P_{k-1}^{i_{1},...,i_{k-2},i_{k-1}}$ and 
$P_{k-1}^{i_{1},...,i_{k-2},i_{k}}$ of the $(d-1)$-dimensional 
projective space $\mathcal{H}_{k-1}^{i_{1},...,i_{k-2}}$.  Define 
$P_{k}^{i_{1},...,i_{k}}$ to be the point where the strict transform
of the line $\mathcal{L}^{i_{1},...,i_{k-2},\{i_{k-1},i_{k}\}}$ meets the
exceptional divisor $\mathcal{H}_{k}^{i_{1},...,i_{k-1}}$.

Denote $\overline{\mathcal{R}}:= \mathcal{R}_{d}$ and
$\eta:=\pi_{1}\circ\pi_{2}\circ \cdots \circ \pi_{d}: \mathcal{R}_{d}
\rightarrow \mathcal{R}$. Let:
$$\overline{\mathcal{H}}:=\eta ^{-1}(P_{0})=
  \bigcup_{\stackrel{k,i_{1},...,i_{k-1}\in \{1,...,d\}}{\mid
      \{i_{1},...,i_{k-1}\} \mid
      =k-1}}\mathcal{H}_{k,d}^{i_{1},...,i_{k-1}}$$
be the exceptional divisor of $\eta$. 
It is by construction a divisor with normal
crossings. The minimal dimension of the centers (see definition
\ref{centre}) of
$\overline{\mathcal{H}}$ 
is $0$ and is attained precisely at the points
$P_{d}^{i_{1},...,i_{d}}$, where  
$(i_{1},...,i_{d})$ varies among the $d!$ permutations of
$(1,...,d)$. Denote this set by $\overline{\mathcal{P}}$.
\medskip

The morphism $\eta$ is isomorphic with a toric morphism obtained from
$\mathcal{Z}(\tilde{W}, \sigma_{0})$ by  a sequence of
star-subdivisions with respect to cones of dimension $d$. One begins
by taking the star-subdivision $\sigma_{0}^{*}$ of $\sigma_{0}$ with
respect to itself. This gives the morphism $\pi_{1}$. Then one
star-subdivides each cone $\sigma_{0}^{i}, \: 1 \leq i \leq d$ (see the
notations preceding before definition \ref{star}),  
and gets $\pi_{2}$. At each one of the next steps, one subdivides
only part of the cones of the maximal dimension. The point is that at
each step it is possible to number canonically the edges of each cone
from 1 to $d$ (the numbering depends on the cone looked upon). Then
one subdivides only cones of the form $(\cdots
((\sigma_{0}^{i_{1}})^{i_{2}})\cdots )^{i_{k}}$, with
$i_{1},...,i_{k}$ pairwise distinct.

If $\tilde{\Sigma}$ denotes the fan obtained by composing all these
subdivisions, the morphism $\eta$ is isomorphic with the toric
morphism $\mathcal{Z}(\tilde{W}, \tilde{\Sigma})\rightarrow
\mathcal{Z}(\tilde{W}, \sigma_{0})$ and $\overline{\mathcal{P}}$ is the
  union of the orbits of dimension 0 of the action of $T_{\tilde{W}}$
  on $\mathcal{Z}(\tilde{W}, \tilde{\Sigma})$.
\medskip

\textbf{Consider then the case in which $c=d$ and $0 < s \leq d-3$.}

The
codimension $1$ part of $(\nu \circ
\mu)^{-1}(\mathrm{Sing}(\mathcal{S}))$ is, by formula
(\ref{preimsing}), the space $\bigcup_{1 \leq i \leq
  s}\mathcal{U}_{i}$. Its center
at the point $P_{0}$ is the space $\bigcap_{1 \leq i \leq
  s}\mathcal{U}_{i}= \mathcal{U}_{[[1,s]]}$, which in the toric
setting is the subspace of the coordinates
$\tilde{U}_{s+1},...,\tilde{U}_{d}$.  The 
intersections $\mathcal{U}_{\{j,l\}} \cap \mathcal{U}_{[[1,s]]}= 
\mathcal{U}_{[[1,s]]\cup \{j,l\}}$ are analytically determined, for
any $j \neq l, \: j,l \in [[s+1, c]]$ and allow to repeat the
construction done in the previous case, looking inside
$\mathcal{U}_{[[1,s]]}$. For example, the axis $\mathcal{L}^{j}$ of
the coordinate $\tilde{U}_{j}$, for $s+1 \leq j \leq d$, can be obtained as:
$\mathcal{L}^{j} = \bigcap_{\stackrel{I \subset \{s+1,...,d\}}
  {\mid I \mid =2, \: j \notin I}} \mathcal{U}_{[[1,s]] \cup I}$. Then,
  instead of blowing-up points, one blows-up smooth varieties of
  dimension $s$, starting with $\mathcal{U}_{[[s+1,d]]}= 
  \bigcap_{\stackrel{ j,l \in [[s+1,d]]}{j \neq l}}
      \mathcal{U}_{\{j,l\}}$.

After $(d-s)$ steps $\pi_{1},..., \pi_{d-s}$, one gets a morphism 
$\eta:= \pi_{1} \circ \cdots \circ \pi_{d-s}:
\mathcal{R}_{d-s}\rightarrow \mathcal{R}$. Denote 
$\overline{\mathcal{R}}:= \mathcal{R}_{d-s}$.  Then  
$\eta ^{-1}(\mathcal{U}_{[[s+1,d]]})$ is the 
exceptional divisor of $\eta$. It is again by construction a divisor
with normal
crossings. The minimal dimension of its centers 
is $s$. Among the points at which the dimension of the center is $s$,
a discrete set $\overline{\mathcal{P}}$ of $(d-s)!$ elements is analytically
distinguished, as formed by the points which are moreover situated on the fiber
$\eta ^{-1}(P_{0})$. Let $\overline{\mathcal{H}}:= \eta
^{-1}(\mathcal{U}_{[[s+1,d]]} \cup \bigcup_{1 \leq i \leq s}
\mathcal{U}_{i})$. Then $\overline{\mathcal{H}}$ is also a divisor
with normal crossings and at the points of $\overline{\mathcal{P}}$
the dimension of its center is $0$. 

In this case, the morphism $\eta$ is isomorphic to a toric morphism
$\mathcal{Z}(\tilde{W}, \tilde{\Sigma})\rightarrow
\mathcal{Z}(\tilde{W}, \sigma_{0})$ obtained by a sequence of
star-subdivisions with respect to cones of dimension $(d-s)$.

\medskip

\textbf{Consider finally the case in which $c < d$.}

The discussion of this case is comparable with the one done in the
previous case. Indeed, one has simply to add $(d-c)$ to the dimensions
of all the varieties which are blown-up. Using the fibers of the
normalization morphism  $\nu: \overline{\mathcal{S}} \rightarrow \mathcal{S}$,
one has at its disposal at each step of blowing-up, canonical families
of projective spaces in which to join points. This gives the
canonical smooth varieties whose strict transforms meet the
exceptional divisors at the new centers of blowing-up. 

At the end, one obtains again a map $\eta: \overline{\mathcal{R}}
\rightarrow \mathcal{R}$ having as exceptional divisor
$\overline{\mathcal{H}}:= \eta
^{-1}(\mathcal{U}_{[[s+1,c]]})\hookrightarrow
  \overline{\mathcal{R}}$. The minimal dimension of its centers  
is now $(s+d-c)$. Among the points at which the dimension of
  the center is $(s+d-c)$, a discrete set $\overline{\mathcal{P}}$ of 
$(c-s)!$ elements is analytically
determined, as formed by the points which are moreover situated on the fiber
$\eta ^{-1}(P_{0})$. We define $\overline{\mathcal{H}}:= \eta
^{-1}(\mathcal{U}_{[[s+1,c]]} \cup \bigcup_{1 \leq i \leq
  s}\mathcal{U}_{i})$.  Then
$\overline{\mathcal{H}}$ is again a divisor with normal crossings and
at the points of $\overline{\mathcal{P}}$, the dimension of its center
is $d-c$. Now $\eta$ is siomorphic to a toric morphism obtained by a
sequence of star-subdivisions with respect to cones of dimension
$(c-s)$. 

\medskip

2) \textbf{Suppose that $s=c-1$}. Then we define $\eta:
  \overline{\mathcal{R}} \rightarrow \mathcal{R}$ to be the blow-up of 
  $\mathcal{U}_{[[1,c]]}$ and $\overline{\mathcal{H}}:=
  \eta ^{-1}(\bigcup_{1 \leq i \leq s}\mathcal{U}_{i})$. The strict
  transform of $\mathcal{U}_{[[1,s]]}=\bigcap_{1 \leq i \leq
  s}\mathcal{U}_{i}$ cuts the fiber $\eta ^{-1}(P_{0})$ in a unique
  point $P_{1}$. We define $\overline{\mathcal{P}}:=\{P_{1}\}$. 

\medskip

3) \textbf{Suppose that $s\in \{c, c-2\}$}.  
Then we do not modify $\mathcal{R}$. So,
$\eta=\mbox{id}_{\mathcal{R}}$ and
$\overline{\mathcal{R}}=\mathcal{R}$. We define then
$\overline{\mathcal{H}}:=\bigcup_{1 \leq i \leq s}\mathcal{U}_{i}$ and 
$\overline{\mathcal{P}}:= \{P_{0}\}$.

\medskip

As all the previous constructions are analytically canonical
presentations of toric constructions, we have obtained the following 
proposition:

\begin{prop} \label{discan}
If $\overline{\mathcal{R}}, \eta, \overline{\mathcal{H}},
\overline{\mathcal{P}}$ are constructed as before, there is a fan
$\tilde{\Sigma}$ obtained from $\sigma_{0}$ by a sequence of
star-subdivisions with respect to cones of dimension $d-c'+s$, such
that $\eta$ is isomorphic to the restriction of the toric map
$\eta_{T}: \mathcal{Z}(\tilde{\mathcal{W}}, \tilde{\Sigma})\rightarrow
\mathcal{Z}(\tilde{\mathcal{W}}, \sigma_{0})$ to a neighborhood of
$\eta_{T} ^{-1}(0)$. Moreover, the components of
$\overline{\mathcal{H}}$ correspond by this isomorphism to orbit
closures of codimension 1, the points of $\overline{\mathcal{P}}$ to
orbits of dimension 0 and at each point of
$\overline{\mathcal{P}}$, the hypersurface $\overline{\mathcal{H}}$
has exactly $c'$ components.
\end{prop}

\section{Reduced Newton polyhedra} \label{NP}

The notions developed in this sections are needed in the next one in
order to define the reduced semigroup of an irreducible quasi-ordinary
polynomial.  

Let $\check{X}$ denote a subset with
$\check{d}$ elements of the variables $X$.

If $\eta \in
\widetilde{\mathbf{C}\{ X \}}$ can be written $ \eta=
(\check{X})^{\check{m}}u(X)$, 
with $ \check{m}\in \mathbf{Q}_{+}^{\check{d}}$, and $u \in
\widetilde{\mathbf{C}\{ 
  X \}}$ verifying $u(\check{X},0)\neq 0$, we say that
$\eta$ \textit{has a dominating exponent with respect to} $\check{X}$,
exponent denoted by $v_{\check{X}}(\eta):= \check{m}$. Here
$u(\check{X},0)$ denotes the series obtained by anulating all the
variables $X_{i}$ which are not in  $\check{X}$.

In fact one can also  extend to this setting the notion of Newton
polyhedron. If 
$\eta \in \widetilde{\mathbf{C}\{ X \}}$, we
define the \textit{reduced Newton polyhedron}
$\mathcal{N}_{\check{X}}(\eta)$
of $\eta$ with respect to $\check{X}$, to be the convex hull in
$\mathbf{R}^{\check{d}}$ of the set 
$\mbox{Supp}_{\check{X}}(\eta)+\mathbf{R}_{+}^{\check{d}}$, where
$\mbox{Supp}_{\check{X}}(\eta)$ denotes the support of $\eta$ written as a
series in the variables $\check{X}$, with coefficients in the
algebra 
$\widetilde{\mathbf{C}\{ X- \check{X}\}}$. 

One can obtain the reduced Newton polyhedron
$\mathcal{N}_{\check{X}}(\eta)$ from the knowledge of the usual Newton
polyhedron $\mathcal{N}_{X}(\eta)$. In order to see it, let
$M=\mathbf{Z}^{d}$ be the lattice of exponents  of the monomials in
$\mathbf{C}\{X\}$ and let  $\check{M} \simeq \mathbf{Z}^{\check{d}}$
be the sublattice of exponents of the monomials in
$\mathbf{C}\{\check{X}\}$. Denote by $\check{p}: M \rightarrow
\check{M}$ the canonical projection of $M$ on $\check{M}$. 

\begin{lemm} \label{projpol}
 One has the equality: $\mathcal{N}_{\check{X}}(\eta)=
 \check{p}_{\mathbf{R}}(\mathcal{N}_{X}(\eta))$. 
\end{lemm}

\textbf{Proof:}  We can suppose without loss of generality that
$\check{X}$ is composed of the first $\check{d}$ variables
$X_{1},...,X_{\check{d}}$. 

Let $\check{A}\in \check{M}_{\mathbf{R}}$ be a vertex of
$\mathcal{N}_{\check{X}}(\eta)$. Then there is a monomial of $\eta$ of
the form $\check{X}^{\check{A}}X_{\check{d}+1}^{j_{1}}\cdots
X_{d}^{j_{d-\check{d}}}$, and so $(\check{A}, j_{1},...,
j_{d-\check{d}})\in \mathcal{N}_{X}(\eta)\Rightarrow \check{A}\in
\check{p}_{\mathbf{R}}(\mathcal{N}_{X}(\eta))$. This implies that 
$ \mathcal{N}_{\check{X}}(\eta) \subset
\check{p}_{\mathbf{R}}(\mathcal{N}_{X}(\eta))$. 

Take now a vertex $\check{A}$ of
$\check{p}_{\mathbf{R}}(\mathcal{N}_{X}(\eta))$. It is immediate to see that it
is the projection of a vertex of $\mathcal{N}_{X}(\eta)$. Choose 
one such vertex by $(\check{A}, j_{1},...,j_{d-\check{d}})\in
M_{\mathbf{Q}}$. Then $(\check{A}, j_{1},...,j_{d-\check{d}})\in
\mbox{Supp}_{X}(\eta)\Rightarrow \check{A}\in
\mbox{Supp}_{\check{X}}(\eta)\Rightarrow  \check{A}\in
\mathcal{N}_{\check{X}}(\eta)$. This implies the reverse inclusion 
$\check{p}_{\mathbf{R}}(\mathcal{N}_{X}(\eta))
\subset\mathcal{N}_{\check{X}}(\eta)$. 
 \hfill $\Box$
\medskip

Even if we do not need it later in this work, we indicate here a
generalization of the dominating exponent, as we think that it has
independent interest. This generalization gives an \textit{intrinsic meaning}
to the notion of reduced Newton polyhedron.

Consider a germ $(\mathcal{H},P)$ of hypersurface with normal
crossings on a smooth variety $\mathcal{V}$. We suppose that the $r$
components $\mathcal{H}_{1},..., \mathcal{H}_{r}$ of $\mathcal{H}$ at
$P$ are taken in a fixed order. 

\begin{prop} \label{invpol}
Let $x=(x_{1},...,x_{n})$ be a system of local coordinates of
$\mathcal{V}$ at $P$, such that
$\mathcal{Z}(x_{i})_{P}=(\mathcal{H}_{i})_{P}, \: 
\forall i \in \{1,...,r\}$.  Denote
$\check{x}:=(x_{1},...,x_{r})$. If $h \in
\mathcal{O}_{\mathcal{V},P}$, then its reduced Newton polyhedron
$\mathcal{N}_{\check{x}}(h)$ depends only on the pair $(h, \mathcal{H})$ and
of the ordering of the components of $\mathcal{H}$ at $P$. 
\end{prop}

\textbf{Proof:} Let $x=(x_{1},...,x_{n})$ and $y=(y_{1},...,y_{n})$ be two
coordinate systems of $\mathcal{V}$, adapted to $\mathcal{H}$ at
$P$. Then, $\forall k \in \{1,...,r\}, \: x_{k}=y_{k}u_{k}$, with
$u_{k}\in \mathcal{O}^{*}_{\mathcal{V},P}$. If $h \in
\mathcal{O}_{\mathcal{V},P}$ and
$h=\sum_{\mbox{Supp}_{\check{x}}(h)}c_{i_{1}\cdots i_{r}}
x_{1}^{i_{1}}\cdots x_{r}^{i_{r}}$ is the expression of $h$ as a
series in the variables $x_{1},...,x_{r}$ with coefficients in
$\mathbf{C}\{ x_{r+1},...,x_{n}\}$, one deduces that $h=
\sum_{\mbox{Supp}_{\check{x}}(h)}(c_{i_{1}\cdots i_{r}}u_{1}^{i_{1}}\cdots
u_{r}^{i_{r}}) y_{1}^{i_{1}}\cdots y_{r}^{i_{r}}$. 

Notice that $c_{i_{1}\cdots i_{r}}u_{1}^{i_{1}}\cdots
u_{r}^{i_{r}}$, seen as an element of $\mathbf{C}\{y_{1},...,y_{n}\}$,  is not
necessarily an element of $\mathbf{C}\{y_{r+1},...,y_{n}\}$. 

Let $m \in \mathbf{N}^{r}$ be a vertex of
$\mathcal{N}_{\check{y}}(h)$ (where
$\check{y}:=(y_{1},...,y_{r})$). It is the exponent of a
$(y_{1},...,y_{r})$-monomial appearing in the development of one of
the terms $(c_{i_{1}\cdots i_{r}}u_{1}^{i_{1}}\cdots
u_{r}^{i_{r}})  y_{1}^{i_{1}}\cdots y_{r}^{i_{r}}$
introduced before. So, for this exponent $(i_{1},...,i_{r})$, we have $m\in
\{(i_{1},...,i_{r})\}+\mathbf{R}_{+}^{r}$. But $(i_{1},...,i_{r})\in
\mbox{Supp}_{\check{x}}(h)\Rightarrow m \in
\mathcal{N}_{\check{x}}(h)$. As this is true for any vertex $m$ of
$\mathcal{N}_{\check{y}}(h)$, it implies that $\mathcal{N}_{\check{y}}(h)
\subset \mathcal{N}_{\check{x}}(h)$. Permuting now the roles of
$x$ and $y$, we get the desired equality
$\mathcal{N}_{\check{x}}(h)= \mathcal{N}_{\check{y}}(h)$. 
\hfill $\Box$
\medskip

We call the invariant Newton polyhedron of the previous proposition
\textit{the Newton polyhedron of $h$ with respect to
  $(\mathcal{H},P)$} and we denote it by
$\mathcal{N}_{(\mathcal{H},P)}(h)$. The function $h$ has a d.e. with
respect to $(\mathcal{H},P)$ (see definition \ref{domexp}) 
if and only if this polyhedron has only
one vertex, which is then equal to $v_{(\mathcal{H},P)}(h)$. So, 
this notion generalizes that of dominating exponent.

\section{The reduced semigroup} \label{redsg}

Recall that the reduced equisingular dimension $c'$ of $\mathcal{S}$ was
defined in section \ref{canoblo}.

If $E$ is a set and $V \in E ^{d}$, we denote by $V'$ the $c'$-tuple
of the first $c'$-coordinates of $V$. Now we particularise the 
constructions of section \ref{NP}  to the case where $\check{X}=X'$. 

The following lemma is an immediate consequence of theorem
\ref{strucsing}. We use the notation $(a)^{j}:=(a,...,a)\in
\mathbf{R}^{j}$. 

\begin{lemm} \label{back}
  For $i \in \{1,...,G-1\}$, one has $A_{i}=A_{i}' \oplus (0)^{d-c'}$,
  and:
  $$A_{G}=  \begin{cases}
              A_{G}' \oplus (0)^{d-c'},  \mbox{ if } s\neq c-2\\
              A_{G}' \oplus (\frac{1}{N_{G}})^{2}\oplus (0)^{d-c}, 
              \mbox{ if } s = c-2
            \end{cases} .$$
  The relations (\ref{relrec}) remain true if one considers the
  vectors $A_{i}'$ and $\overline{A}_{i}'$ instead of the vectors  
  $A_{i}$ and $\overline{A}_{i}$.
\end{lemm}

We imitate now the construction of the semigroups
$\Gamma_{\mathcal{N}}(f)$ and 
$\Gamma_{\mathcal{D}}(f)$, by introducing the following subsets of
$\mathbf{Q}^{c'}_{+}$:
$$\Gamma'_{\mathcal{D}}(f):= \{ v_{X'}(h(\xi)), h \in
\mathbf{C}\{X\}[Y]-(f), \: h(\xi)  
   \: \mbox{ has a d.e. with respect to } X'\},$$
$$ \Gamma'_{\mathcal{N}}(f):= \{ A' \in  \mathbf{Q}_{+}^{c'}, \: 
    A'\mbox{ is a vertex of } \mathcal{N}_{X'}(h(\xi)), \: 
    h \in \mathbf{C}\{X\}[Y]-(f) \}.$$ For the same reasons as before,
    they are additive semigroups. 

We have the following proposition, to be compared with propositions
\ref{gensgr} and \ref{samesgr}:

\begin{prop} \label{samered}
One has the equality of semigroups:
  $$\Gamma'_{\mathcal{N}}(f)=\Gamma'_{\mathcal{D}}(f)=
  \mathbf{N}^{c'}+ \mathbf{N} \overline{A}'_{1} 
  + \cdots + \mathbf{N} \overline{A}'_{G}.$$
\end{prop}

\textbf{Proof:} Take any $A \in \mathbf{N}^{c'}$ and
$(i_{0},...,i_{G-1})\in \mathbf{N}^{G}$. As \linebreak 
$v_{X'}((X')^{A}(f_{0}(\xi))^{i_{0}}\cdots (f_{G-1}(\xi))^{i_{G-1}})=
A+ i_{0}\overline{A}_{1}'+\cdots + i_{G-1}\overline{A}_{G}'$, one gets
the inclusions 
$\mathbf{N}^{c'}+
\mathbf{N}\overline{A}_{1}'+\cdots+ \mathbf{N}\overline{A}_{G}' \subset
\Gamma_{\mathcal{D}}'(f)\subset \Gamma_{\mathcal{N}}'(f)$. 

Suppose now that $A'\in \Gamma_{\mathcal{N}}'(f)$. Then there is an
element $h \in \mathbf{C}\{X\}[Y]-(f)$ such that $A'$ is a vertex of
$\mathcal{N}_{X'}(h(\xi))$. By lemma \ref{projpol}, there is a vertex
of $\mathcal{N}_{X}(h(\xi))$ of the form
$(A',j_{1},...,j_{d-c'})$. But, by proposition \ref{gensgr}, this
vertex is an element of $\mathbf{N}^{d}+
\mathbf{N}\overline{A}_{1}+\cdots+ \mathbf{N}\overline{A}_{G}$, which
implies that $A'$ is an element of $\mathbf{N}^{c'}+
\mathbf{N}\overline{A}_{1}'+\cdots+ \mathbf{N}\overline{A}_{G}'$.
\hfill $\Box$
\medskip

This motivates the following definition:

\begin{defin} \label{sgrextr}
The semigroup $\Gamma'_{\mathcal{N}}(f)=\Gamma'_{\mathcal{D}}(f)$ is
called \textbf{the reduced semigroup of} $\mathcal{A}$ 
\textbf{with respect to} $f$, denoted $\Gamma'(f)$.
\end{defin}

If $(\Gamma,+)$ is a finitely generated abelian semigroup without
torsion, let  $\hat{\Gamma}$ be the lattice generated by $\Gamma$. 
Denote by $\sigma(\Gamma)$ the 
convex cone generated by $\Gamma$ in $\hat{\Gamma}_{\mathbf{R}}$.

Particularize this to the reduced semigroup of $f$. As
$\Gamma'(f)\subset \mathbf{Q}_{+}^{c'}$, the cone $\sigma(\Gamma'(f))$
is strictly convex. Let $u^{1},...,u^{c'}$  be the 
smallest non-zero elements of $\Gamma'(f)$ situated 
on the edges of $\sigma(\Gamma'(f))$. The following lemma shows that
almost all the vectors $\overline{A}_{1}',...,\overline{A}_{G}'$ of
$\mathbf{Q}_{+}^{c'}$ are determined by the isomorphism type of $\Gamma'(f)$:

\begin{lemm} \label{unigen}
For any $j \geq 1$, if $\alpha_{1},...,\alpha_{j-1}$ are already defined 
and verify $\Gamma'(f) \neq \mathbf{N}u^{1}+\cdots +
\mathbf{N}u^{c'}+ 
\mathbf{N}\alpha_{1}+\cdots + \mathbf{N}\alpha_{j-1}$, 
there exists a unique smallest element $\alpha_{j}$ of $\Gamma'(f)$ 
not contained in the semigroup $\mathbf{N}u^{1}+\cdots +\mathbf{N}u^{c'}+ 
\mathbf{N}\alpha_{1}+\cdots + \mathbf{N}\alpha_{j-1}$. Define 
$g \in \mathbf{N}$ by $\Gamma'(f) = \mathbf{N}u^{1}+\cdots +
\mathbf{N}u^{c'}+ 
\mathbf{N}\alpha_{1}+\cdots + \mathbf{N}\alpha_{g}$ and $\epsilon \in
\{ 0,1\}$ to be either $0$ if $s \neq c-2$ or $1$ if $s=c-2$. Then $g
\in \{G-\epsilon,G\}$ and, after possibly 
permuting $u^{1},...,u^{c'}$, the components of
$\alpha_{1},...,\alpha_{G-\epsilon}$ 
written in the basis $u^{1},...,u^{c'}$ coincide with the vectors  
$\overline{A}_{1}',...,\overline{A}_{G-\epsilon}'$ of
$\mathbf{Q}_{+}^{c'}$ associated  
to any normalized qo-defining polynomial of $\mathcal{S}$.
\end{lemm}

\textbf{Proof:} By lemma \ref{back}, $\forall k \in
\{1,...,G-\epsilon\}$, one has $N_{k}= \mbox{min}\{j \in
\mathbf{N}^{*}, \: j \overline{A}_{k}'\in \mathbf{Z}^{c'}+ 
\mathbf{Z}\overline{A}_{1}'+\cdots+
\mathbf{Z}\overline{A}_{k-1}'\}$. As $N_{k}>1$, we get $
\overline{A}_{k}'\notin  \mathbf{N}^{c'}+ \mathbf{N}\overline{A}_{1}'+\cdots+
\mathbf{N}\overline{A}_{k-1}'$. But if $\overline{A}'\in \Gamma'(f)$
and $\overline{A}'\notin \mathbf{N}^{c'}+ \mathbf{N}\overline{A}_{1}'+\cdots+
\mathbf{N}\overline{A}_{k-1}'$, then $\overline{A}'=A_{0}'+
\sum_{j=1}^{G}i_{j}\overline{A}_{j}'$, with $A_{0}'\in \mathbf{N}^{c'}$
and at least one of $i_{k},...,i_{G}$ is non-zero. This implies that
$\overline{A}' \geq \overline{A}_{k}'$, and so $\overline{A}_{k}'$ is
the unique smallest element of $\Gamma'(f)$ not contained in the
semigroup $\mathbf{N}^{c'}+ \mathbf{N}\overline{A}_{1}'+\cdots+
\mathbf{N}\overline{A}_{k-1}'$. As $(u ^{1},...,u ^{c'})$ form
obviously a permutation of the canonical generators of the semigroup
$(\mathbf{N}^{c'},+)$, this proves the lemma. 
\hfill $\Box$
\medskip

This lemma will be used in the passage from the analytical invariance
of the reduced semigroup to the analytical invariance of the
normalized characteristic exponents (Corollary \ref{invchar}).

\section{The main results} \label{results}

In the sequel we suppose that the germ $(\mathcal{S},0)$ is irreducible, 
quasi-ordinary of 
dimension $d \geq 2$ and embedding dimension $d+1$. Moreover, we suppose that 
$0$ is not smooth on $\mathcal{S}$. Let $f$ be a qo-defining polynomial of 
$\mathcal{S}$. With our hypothesis, $f$ is irreducible.

Define 
$\nu:\overline{\mathcal{S}}\longrightarrow \mathcal{S}$ to be 
the normalization morphism of $\mathcal{S}$, studied in section
\ref{norgerm}.  
Let 
$\mu: \mathcal{R}\longrightarrow \overline{\mathcal{S}}$ be 
the orbifold map of $\overline{\mathcal{S}}$, introduced in 
section \ref{orbmap}. 
Let also 
$\eta:\overline{\mathcal{R}}\longrightarrow \mathcal{R}$ be the 
canonical modification of $\mathcal{R}$ defined in section 
\ref{canoblo}.

Denote $\theta := \nu \circ \mu \circ \eta: \overline{\mathcal{R}}
\rightarrow \mathcal{S}$. Let $\theta^{*}$ be the corresponding morphism 
of sheaves of local algebras. By construction, the morphism $\theta$ depends
only on the analytical  
type of $\mathcal{S}$ and not on any particular qo-defining polynomial.
Comparing the construction done in the previous section and the
definition of $\mathcal{Z}'$ given in section \ref{canoblo}, we get:
$$\overline{\mathcal{H}}= \theta ^{-1}(\mathcal{Z}').$$

Let $P$ be a 
point of $\overline{\mathcal{P}}$. By proposition \ref{discan}, the
hypersurface $ \overline{\mathcal{H}}$ has $c'$ components at $P$.
The following definition generalizes the first one presented in 
section 1 for plane curves:

\begin{defin} \label{semigrP}
The \textbf{reduced semigroup of $\mathcal{S}$ 
with respect to $P$},
 denoted 
 $\Gamma_{P}'(\mathcal{S})$, is the following subsemigroup of 
$(\mathbf{N}^{c'},+)$:
$\Gamma_{P}'(\mathcal{S}):=
\Gamma_{\overline{\mathcal{H}},P}
(\theta^{*}(\mathcal{A})_{P}).$
\end{defin}

The semigroup $\Gamma_{P}'(\mathcal{S})$ depends obviously only on the
choice of the point $P$. Theorem \ref{isosgsg} shows that in fact,
up to isomorphism, this semigroup
is independent of $P$, and so  it is 
an \textit{analytical invariant} of $\mathcal{S}$: 

\begin{theor} \label{isosgsg}
Let $f$ be a quasi-ordinary defining polynomial of $\mathcal{S}$. For  
every  point $P\in \overline{\mathcal{P}}$,  
the image of the composition 
$\mathcal{E}'(f)  \rightarrow \mathcal{A}  \rightarrow \mathcal{O}_{
\overline{\mathcal{R}},P}$ of the restriction mapping and of 
 $\theta^{*}$
is contained in the set
$\mathcal{E}_{\overline{\mathcal{H}},P}(\theta^{*}(\mathcal{A})_{P})$.  
It induces a 
well-defined mapping $\Phi_{P}$ which realises an isomorphism of semigroups:
$$ \begin{array}{lll}\Phi_{P}: &
              \Gamma'(f) & \longrightarrow \Gamma'_{P}(\mathcal{S})\\
              & v_{X'}(h(\xi)) & \longrightarrow
              v_{\overline{\mathcal{H}},P}(\theta^{*}
                 (h\mid_{\mathcal{S}})_{P})
            \end{array}.$$
Here $h$ varies through $\mathcal{E}'(f)$.
\end{theor}

We prove this theorem in the next section. It obviously generalizes
the theorem \ref{mths}.

As the left-hand semigroup does not 
depend on the choice of the point $P\in \overline{\mathcal{P}}$,  
and the right-hand one does not depend on the 
choice of qo-defining polynomial $f$, we get:

\begin{corollary} \label{invsem}
As an abstract semigroup, 
$\Gamma'(f)$ does not depend on the chosen defining polynomial $f$ of 
 $\mathcal{S}$.  We call it \textbf{the reduced
semigroup of } $\mathcal{S}$, denoted  $\Gamma'(\mathcal{S}).$
\end{corollary}

As a by-product of the proof of theorem \ref{isosgsg}, we get a way to 
associate to some elements of $\mathcal{A}$ a value in the semigroup 
$\Gamma'(\mathcal{S})$:

\begin{corollary} \label{imsg}
Let $f$ be a qo-defining polynomial of the germ  
$\mathcal{S}$ and $\xi \in R(f)$. If $h \in \mathcal{E}'(f)$, 
then the dominating exponent $v_{X'}(h(\xi))$, seen as an element of the  
abstract semigroup $\Gamma'(\mathcal{S})$, depends only on the image 
$h \mid_{\mathcal{S}} \in \mathcal{A}$, and not on the choice of $f$.
\end{corollary}

%If $\mbox{Sing}(\mathcal{S})$ has two irreducible components, then the 
%image of $\mathcal{E}(f)$ in $\mathcal{A}$ obtained by restriction is 
%independent of the defining $f$. This fact is no longer true if 
%$\mbox{Sing}(\mathcal{S})$ is irreducible (look at  
%$f: = Y^2 - X_{1}X_{2}^{3}, \: h:=Y\in \mathcal{E}(f)$ and change the 
%quasi-ordinary projection).

Our main application of theorem \ref{isosgsg} is the following one:

\begin{corollary} \label{invchar}
The characteristic exponents of a normalized 
qo-defining polynomial $f$ of the germ $\mathcal{S}$ 
are analytical invariants of $\mathcal{S}$ (we recall that we suppose
condition (\ref{lexord}) is verified).
\end{corollary}

We postpone the proof to section \ref{pfcor}.
\medskip

\textbf{Remark: } From Corollary \ref{invchar} one can deduce also the
analytical invariance of the semigroup $\Gamma(f)$, as was done in
\cite{GP 00} using the inversion formulae, expressing arbitrary
characteristic exponents in terms of normalized ones.

\section{Proof of the theorem \ref{isosgsg}} \label{pfs}

Recall from section \ref{Strucsing} that $s$ denotes the number of 
components of $\mathrm{Sing}(\mathcal{S})$
which have codimension 1 in $\mathcal{S}$.

Denote $\overline{\psi}:= \psi \circ \theta: \overline{\mathcal{R}} 
\rightarrow \mathbf{C}^{d}$. We have the following commutative
diagram: 

$$\xymatrix{
     P \ar[r] & \overline{\mathcal{P}} \ar[r] &
     \overline{\mathcal{H}}= \theta^{-1}(\mathcal{Z}') \ar[d] & & &  
        \mathcal{Z}' \ar[d] \ar[r] & \mbox{Sing}(\mathcal{S}) \\
    &  & \overline{\mathcal{R}}  \ar@/^2pc/[rrr]^{\theta}
            \ar[r]^{\eta} \ar[drrr]_{\overline{\psi}} & 
    \mathcal{R} \ar[r]^{\mu} \ar[drr] & \overline{\mathcal{S}} \ar[r]^{\nu} 
      \ar[dr] &  \mathcal{S} \ar[d]^{\psi} \\
      &   &  & & & \mathbf{C}^{d} & \ar[l] \mathcal{D}'}$$

As in section \ref{simple}, take $\xi:= Y\mid_{\mathcal{S}}$. One uses
again formula (\ref{relfund}). 

One can choose as representative of $\overline{\psi}$ a localisation to an 
open set of a toric morphism. 
Indeed, following the sections \ref{norgerm}, \ref{orbmap} and 
\ref{canoblo}, one can realize representatives of $\nu, \mu, \eta$ 
as toric morphisms $\nu: \mathcal{Z}(W_{G}, \sigma_{0})\rightarrow
\mathcal{Z}(W, \sigma_{0}), \: \mu: \mathcal{Z}(\tilde{W},
\sigma_{0})\rightarrow 
\mathcal{Z}(W_{G}, \sigma_{0}), \: \eta: \mathcal{Z}(\tilde{W},
\tilde{\sigma})\rightarrow 
\mathcal{Z}(\tilde{W}, \sigma_{0})$. With such representatives of the
morphisms, the  
point $P$ is an orbit of dimension $0$ and $\overline{\mathcal{H}}$ is a union
of closures of the orbits of codimension 1. Let $T:=(T_{1},...,T_{d})$
be the toric coordinates of $\overline{\mathcal{R}}$ centered at
$P$. They are adapted to $\overline{\mathcal{H}}$ at $P$. So,
$\overline{\mathcal{H}}_{P}= Z(T_{1})\cup \cdots \cup Z(T_{c'})$. 
With such coordinates, the morphism $\overline{\psi}^{*}$ is 
monomial. 

Let $\tau_{1},...,\tau_{d} \in \tilde{M}$ be such that $T_{i}=\chi
^{\tau_{i}}, \forall i \in \{1,...,d\}$. Denote by $T'$ the set of
coordinates $T_{1},...,T_{c'}$ and by $p'$ the projection
of $\tilde{M}$ on the sublattice $\tilde{M}'$ generated by
$\tau_{1},..., \tau_{c'}$. If $\phi: \tilde{W}\rightarrow W$ is the
morphism obtained by composing the changes of lattices of sections
\ref{norgerm}, \ref{orbmap}, we see 
that $p'\circ \check{\phi}: M \rightarrow \tilde{M}'$  depends only on
the restriction of $p'\circ \check{\phi}$ to $M'$ and that the
restriction $p'\circ \check{\phi}:M' \rightarrow  \tilde{M}'$ is an
isomorphism. 
Using formula (\ref{relfund}), we see that $\Phi_{P}$ 
\textit{is injective}.

In order to prove its surjectivity, we must show that if 
$h \in \mathbf{C}\{X\}[Y]$ verifies 
$\theta^{*}(h\mid_{\mathcal{S}})_{P} \in 
\mathcal{E}_{\overline{\mathcal{H}},P}(\theta^{*}(\mathcal{A})_{P})$, 
then one can find another element 
$\tilde{h} \in \mathcal{E}'(f)$ such that
$v_{\overline{\mathcal{H}},P}(\theta^{*} 
(\tilde{h} \mid_{\mathcal{S}})_{P})=
v_{\overline{\mathcal{H}},P}
(\theta^{*}
(h\mid_{\mathcal{S}})_{P})$. As $f\mid_{\mathcal{S}}=0$, we can suppose 
that $\mbox{deg}(h)< \mbox{deg}(f)$, after possibly making the euclidian 
division of $h$ by $f$. We consider then a complete system 
$(f_{0},...,f_{G})$ of semiroots of $f$ and the $(f_{0},...,f_{G})$-adic 
expansion of $h$, which by our hypothesis is of the form 
$h= \sum c_{i_{0}\cdots i_{G-1}} f_{0}^{i_{0}}\cdots
f_{G-1}^{i_{G-1}}$. 

Formula (\ref{relfund}) implies $v_{\overline{\mathcal{H}},P}(\theta^{*} 
(\tilde{h} \mid_{\mathcal{S}})_{P})=
v_{T'}(\overline{\psi}^{*}(h(\xi)))$. By lemma \ref{disjpol}, this
vector is the image by $p'$ of a vertex $\mu$ of
$\mathcal{N}_{T}(\overline{\psi}^{*}(h(\xi)))$. By lemma \ref{transf},
$\mu=\check{\phi}(m)$, where $m$ is a vertex of
  $\mathcal{N}_{X}(h(\xi))$. Now, lemma \ref{corlemm} shows that $m$
  is a vertex of one of the polyhedra $\mathcal{N}_{X}(c_{i_{0}\cdots
    i_{G-1}} (f_{0}(\xi))^{i_{0}}\cdots 
(f_{G-1}(\xi))^{i_{G-1}})$. In particular, it is of the form
$v_{X}(X^{A}  (f_{0}(\xi))^{i_{0}}\cdots 
(f_{G-1}(\xi))^{i_{G-1}})$, where $A$ is a vertex of
$\mathcal{N}_{X}(c_{i_{0}\cdots  i_{G-1}})$. So, $\mu=
p'(v_{\overline{\mathcal{H}},P}(\theta ^{*}(X^{A}f_{0}^{i_{0}}\cdots
f_{G-1}^{i_{G-1}}\mid_{\mathcal{S}})_{P})$. 
But 
$X^{A} f_{0}^{i_{0}}\cdots f_{G-1}^{i_{G-1}} \in \mathcal{E}'(f)$,
showing that we can take $\tilde{h}:= X^{A} f_{0}^{i_{0}}\cdots
f_{G-1}^{i_{G-1}}$. This  
proves that $\Phi_{P}$ \textit{is  surjective}. \hfill $\Box$
\medskip

\section{Proof of the corollary \ref{invchar}} \label{pfcor}

We need first some classical results about germs of plane curves and
about the equisingularity of plane sections of a germ of
hypersurface. 

One has the following classical theorem, attributed sometimes to
M.Noether (see proposition 6.5 in \cite{PP 99}):

\begin{prop} \label{intcoinc} Let $x,y$ be two indeterminates. 
Let $f,g \in \mathbf{C}[[x]][y]$ be irreducible unitary
polynomials. Denote by $(f,g)$ their intersection number and by:
$$K(f,g):=\mathrm{ max } \{v_{x}(\xi-\eta), \: \xi \in R(f), \: \eta \in
R(g)\}$$ the \textbf{exponent of coincidence} of $f$ and $g$, where
their roots are seen as Newton-Puiseux series. If $A_{1},...,A_{G}$
are the characteristic exponents of $f$, that $A_{G+1}:=+\infty$ and
$k \in \{0,...,G\}$ is the smallest integer so that $K(f,g)<A_{k+1}$,
then:
$$\frac{(f,g)}{d_{Y}(f)d_{Y}(g)}=\frac{\overline{A}_{k}}{N_{1}\cdots N_{k-1}}+ 
  \frac{K(f,g)-A_{k}}{N_{1}\cdots N_{k}}.$$
\end{prop}

The following proposition is an easy consequence of the also classical
\textit{inversion formulas} for plane germs (see proposition 4.3 in
\cite{PP 99}):

\begin{prop} \label{1exp}
Let $C$ be a non-regular germ of irreducible plane curve and $a_{1}>1$
be its first characteristic exponent in generic coordinates. Then the
possible values of the first characteristic exponent in various
coordinate systems span the set $\{a_{1}, \frac{1}{a_{1}}\}\cup
\{\frac{1}{m}, \: m \leq [a_{1}]\}$. Moreover, if the
first characteristic exponents of $C$ in two coordinate systems
coincide, then the whole sequence of such exponents coincide.
\end{prop}

The following proposition is also classical (see \cite{Z 70}, \cite{Z
  86}):

\begin{prop} \label{equicurv}
Let $C$ and $C'$ be two germs of reduced plane curves. The following
conditions are equivalent:

1) There is a 1-to-1 correspondence $C_{i}\leftrightarrow C'_{i}$
   between the components of $C$ and $C'$  such that $\forall i$,
   $C_{i}$ and $ C'_{i}$ have the same sequence of characteristic
   Newton-Puiseux exponents and $\forall i\neq j$, one has the
   equality $(C_{i}, C_{j})=(C_{i}', C_{j}')$ of intersection numbers.

2) The germs $C$ and $C'$ have isomorphic processes of embedded
   resolution by blowing-ups. 

3) The dual graphs of their total
   transforms by the minimal embedded resolution morphisms are
   isomorphic.
\end{prop}

This motivates the following definition:

\begin{defin} \label{defequis}
The germs $C$ and $C'$ are said to be \textbf{equisingular} if the
equivalent conditions of the previous proposition are verified.
\end{defin}

Let $(\mathcal{V},0)$ be a germ of irreducible complex analytical space
of dimension $d$ and embedding dimension $d+1$. Let $\mathcal{H}$ be
an irreducible germ of hypersurface on $\mathcal{V}$. Consider an
embedding $E:(\mathcal{V},0)\rightarrow (\mathbf{C}^{d+1},0)$ and
denote by $\psi_{E}$ the projection on the first $d$
coordinates. Suppose that $\psi_{E}$ is finite. Let $P$ be a smooth
point of $\psi_{E}(\mathcal{H})$, and let $L$ be a smooth germ of
curve in $\mathbf{C}^{d}\times\{0\}$, transversal to
$\psi_{E}(\mathcal{H})$ at $P$. Let $Q$ be any point of $\psi_{E}
^{-1}(P)\cap \mathcal{V}$. Zariski \cite{Z 70} proves the following:

\begin{theor} \label{equisect}
The equisingularity type of the germ of plane curve $(\psi_{E}
^{-1}(L)\cap \mathcal{V},Q)$ depends only on
$(\mathcal{V},\mathcal{H},0)$ and not on the choices of $E,P,Q,L$.
\end{theor}

For more informations about the notion of equisingularity, one can
consult Teissier \cite{T 75}.

\medskip

Let us pass now to the proof of the corollary. 

Suppose that $f$ is a \textit{normalized} qo-defining polynomial of
$\mathcal{S}$. Lemma \ref{back} shows that the characteristic
exponents $A_{1},...,A_{G}$ of $f$ are known , once
$A_{1}',...,A_{G}'$ and $N_{G}$ are known. The same lemma shows that
the knowledge of these last quantities is equivalent with the one of
$\overline{A}_{1}',..., \overline{A}_{G-1}', A_{G}'$ and $N_{G}$. 
By lemma \ref{unigen}, the abstract semigroup
$\Gamma'(\mathcal{S})$ allows to determine the vectors
$\overline{A}_{1}',...,\overline{A}_{G-1}'$ of
$\mathbf{Q}_{+}^{c'}$. The same lemma shows that $\forall \: i \in
\{1,...,G-1\}$, the number $N_{i}$ is the index of $\mathbf{Z}^{c'}+
\mathbf{Z} \overline{A}_{1}'+\cdots +\mathbf{Z} \overline{A}_{i-1}'$
in $\mathbf{Z}^{c'}+
\mathbf{Z} \overline{A}_{1}'+\cdots +\mathbf{Z} \overline{A}_{i}'$,
which proves its analytical invariance. Moreover, the same is true for
$\overline{A}_{G}'$ and $N_{G}$, excepted when $s=c-2$. So, the
corollary is proved for $s \neq c-2$.

\textbf{Suppose now that $s=c-2$}. We will show that in addition to  
$\overline{A}_{1}',..., \overline{A}_{G-1}'$, one can also determine
$A_{G}'$ and $N_{G}$ from the analytical type of $(\mathcal{S},0)$. In
order to do it, we shall use no more the analytical invariance of
the semigroup $\Gamma'(f)$. Instead, we shall use the results
presented at the beginning of the section.

In what follows, if $P$ is a point of $\mathbf{C}^{d}$ and $I\subset
[[1,d]]$, we denote by $\mathcal{D}_{I}^{P}$ the $\mid I
\mid$-codimensional affine subspace of $\mathbf{C}^{d}$ passing
through $P$ and parallel to $\mathcal{D}_{I}$.

By theorem \ref{strucsing}, the variety
  $\mathcal{Z}_{\{c-1,c\}}$ is analytically distinguished as the only
  component of $\mathrm{Sing}(\mathcal{S})$ of codimension 2 in
  $\mathcal{S}$. 
We can obtain $N_{G}$ from the analytical
structure of 
$\mathcal{S}$, as \textit{the least number of blowing-ups of the strict
transforms of $\mathcal{Z}_{\{c-1,c\}}$, needed to desingularize
$\mathcal{S}$ at their
  generic points}. This is an immediate consequence of theorem
  \ref{strucsing}, point 3).

\medskip

It remains to show that $A'_{G}$ is analytically determined by
$\mathcal{S}$. We do it by showing that $A ^{i}_{G}\in \mathbf{Q}_{+}$
is analytically determined for any $i \in \{1,...,c-2\}$.

\textbf{1) Consider the case where $c=3$ and $G=1$}. Then
$A_{1}=(A_{1}^{1}, \frac{1}{N_{1}},  \frac{1}{N_{1}})\oplus (0)^{d-3}$, with
$B_{1}^{1}:=N_{1}A_{1}^{1} \in \mathbf{N}^{*}-\{1\}$ (see theorem
\ref{strucsing}).  Take $P \in
\mathcal{D}_{1}- (\mathcal{D}_{2}\cup \mathcal{D}_{3})$. Denote $C^{P}:= \psi
^{-1}(\mathcal{D}_{[[2,d]]}^{P})\cap \mathcal{S}$. The intersection
$\psi ^{-1}(P)\cap \mathcal{S}$  is reduced to one point  $Q$. By
theorem \ref{equisect}, the equisingularity type of $(C^{P},Q)$ is an
analytical invariant of $(\mathcal{S},0)$. 
The germ $(C^{P},Q)$ has $\mbox{gcd}(N_{1}, B_{1}^{1})$ irreducible
components, and so this number is also analytically determined by
$(\mathcal{S},0)$.  

\textbf{If $\mbox{gcd}(N_{1},B^{1}_{1})=1$}, then
$\frac{B^{1}_{1}}{N_{1}}$ is the 
unique Newton-Puiseux characteristic exponent of $C^{P}
\hookrightarrow \psi ^{-1} (\mathcal{D}_{[[2,d]]}^{P})$ with respect
to the coordinates $(X_{1},Y)$ of the plane
$\psi ^{-1} (\mathcal{D}_{[[2,d]]}^{P})$. As $A_{1}^{1}\neq
\frac{1}{N_{1}}$, neither $A_{1}^{1}$ 
nor $\frac{1}{A_{1}^{1}}$ is an integer. By proposition \ref{1exp}, we
deduce that $A_{1}^{1}$ is determined as the unique first
characteristic exponent of $(C^{P},Q)$, for varying coordinate
systems, having $N_{1}$ as the denominator of its irreducible form.

\textbf{If $\mbox{gcd}(N_{1},B^{1}_{1})>1$}, then take two irreducible
components of $C^{P}$ at $Q$. By proposition \ref{intcoinc}, 
their intersection number at $Q$ is
$\frac{N_{1}}{\mbox{gcd}(N_{1}, B^{1}_{1})}\cdot
\frac{B^{1}_{1}}{\mbox{gcd}(N_{1}, 
  B^{1}_{1})}$. As $N_{1}$ and $\mbox{gcd}(N_{1}, B^{1}_{1})$ are
already determined, 
this determines $B^{1}_{1}$, and consequently $A ^{1}=\frac{B^{1}_{1}}{N_{1}}$.
\medskip

\textbf{2) Consider the case where $c=3$ and $G>1$}. Then the
characteristic monomials of $f$ are
$X_{1}^{A_{1}^{1}},...,X_{1}^{A_{G-1}^{1}},
X_{1}^{A_{G}^{1}}X_{2}^{\frac{1}{N_{G}}}X_{3}^{\frac{1}{N_{G}}}$.  We
take the same notations as in the previous case. If $A_{G}^{1}$ is a
characteristic exponent of the components of the germ $(C^{P},Q)$ in
the coordinates $(X_{1}, Y)$, then  proposition \ref{1exp} shows that
it is determined by the 
equisingularity type of $(C^{P},Q)$, as $A_{1}^{1},..., A_{G-1}^{1}$
are already known.  If
$A_{G}^{1}$ is not a characteristic exponent of the components of
$(C^{P},Q)$, then there are at least two such components. We look  at
the intersection number $I$ of any two of them.  
By proposition \ref{intcoinc}, we have $\frac{I}{N_{1}\cdots
  N_{G-1}}= N_{G-1} \overline{A}^{1}_{G-1}- A_{G-1}^{1}+A_{G}^{1}$,
which determines $A_{G}^{1}$ from the knowledge of
$A_{1}^{1},...,A_{G-1}^{1}, I$.  Indeed, lemma \ref{back} shows that
$N_{1},...,N_{G-1}$ can be deduced from $A_{1}^{1},...,A_{G-1}^{1}$.

\medskip

\textbf{3) Consider the case where $c\geq 4$ and $G\geq1$}. Then the
characteristic monomials of $f$ are
$X_{1}^{A_{1}^{1}}\cdots
X_{c-2}^{A_{1}^{c-2}},...,X_{1}^{A_{G-1}^{1}}\cdots X_{c-2}^{A_{G-1}^{c-2}},
X_{1}^{A_{G}^{1}}\cdots
X_{c-2}^{A_{G}^{c-2}}X_{c-1}^{\frac{1}{N_{G}}}X_{c}^{\frac{1}{N_{G}}}$.
We want to show that for all $i \in \{1,...,c-2\}$, $A_{G}^{i}\in
\mathbf{Q}_{+}$ is determined by 
the analytical structure of $(\mathcal{S},0)$.

Let $P(t_{1},...,t_{i-1},0,t_{i+1},...,t_{d})$ be a point of
$\mathcal{D}_{i}-\cup_{j\neq i}\mathcal{D}_{j}$.   Define $C_{i}^{P}:=
\psi ^{-1}(\mathcal{D}_{[[1,d]]-\{ 
  i\}}^{P})\cap \mathcal{S}$.  It is a curve embedded in the plane
$\psi ^{-1}(\mathcal{D}_{[[1,d]]-\{i\}}^{P})$. Localize it at any
point $Q \in \psi ^{-1}(P)\cap \mathcal{S}$ and look at the
equisingularity type of the germ $(C^{P},Q)$. In particular, at the
characteristic exponents of the irreducible components and at the
intersection numbers of the pairs of components, as in the previous two
cases. Each component has the same characteristic exponents as the
curve with Newton-Puiseux series $X_{i}^{A_{1}^{i}}+\cdots +
X_{i}^{A_{G}^{i}}$. These characteristic exponents form a
(possibly strict) subset of $\{A_{1}^{i},..., A_{G}^{i}\}$.  We
subdivide this case in two subcases, analogous with the cases 1),
respectively 2) treated before.

\textbf{Suppose that $A_{1}^{i},..., A_{G-1}^{i}\in \mathbf{N}$}. Then
$B_{G}^{1}:= N_{G}A_{G}^{1}\in \mathbf{N}$ and we determine
$A_{G}^{1}$ as in case 1), using the fact that $s=c-2$ implies $B_{G}^{1}>1$.

\textbf{Suppose that at least one of $A_{1}^{i},..., A_{G-1}^{i}$ is
  in  $\mathbf{Q}-\mathbf{N}$}. Then, as $A_{1}^{i},..., A_{G-1}^{i}$
  are known, we can determine the characteristic exponents among
  them. Then we determine $A_{G}^{1}$ as in case 2).
\hfill $\Box$

\section{Comparison with the 2-dimensional case} \label{comp2}

Let us now compare this work with the paper \cite{PP 02}, in which 
we had obtained the analytical invariance of the semigroup and of the
normalized 
characteristic exponents in the case of quasi-ordinary surfaces.

In \cite{PP 02} the strategy of proof was the same. The morphism $\theta:
\overline{\mathcal{R}} \rightarrow \mathcal{S}$ was obtained also as a
composition $\theta= \nu \circ \mu \circ \eta$, with $\eta:
\overline{\mathcal{S}} \rightarrow \mathcal{S}$ the normalization
morphism of $\mathcal{S}$. The proof of the isomorphism of the
semigroups was basically the same as here.

The main difference is that we defined $\mu$ to be \textit{the minimal
  resolution} of $\overline{\mathcal{S}}$. The last morphism $\eta$
  was either the identity or the blow-up of a point. We took as
  hypersurface $\overline{\mathcal{H}}$ the full preimage $\theta
  ^{-1}(\mbox{Sing}(\mathcal{S}))$ of the singular locus of
  $\mathcal{S}$ and as set $\overline{\mathcal{P}}$ the union of its
  singular points. With the exception of the case when $\mathcal{S}$
  was analytically isomorphic with the germ of a quadratic cone at its
  vertex, we showed that $\overline{\mathcal{P}} \neq
  \emptyset$. At the points of  $\overline{\mathcal{P}}$, the curve
  $\overline{\mathcal{H}}$ had two local components, and so \textit{it was not
  needed to consider reduced semigroups}, as done in the present paper.

The obstruction to extend that method to higher dimensions was that in
general \textit{one has no more unicity of the minimal resolution of
$\overline{\mathcal{S}}$}. Even the existence of some canonical
non-minimal resolution would have been enough, if it could be obtained
by a toric morphism once $\overline{\mathcal{S}}$ was presented as a
germ of toric variety in the way explained in section \ref{norgerm}. 

The solution of our problem of extension to higher dimensions came
once we looked for a morphism $\mu$ \textit{without asking it to be
birational}. Indeed, the attributes of $\mu$ that were important for us
were that its source was smooth and that it could be presented as a
toric morphism once $\overline{\mathcal{S}}$ was presented as a
germ of toric variety in the way explained in section
\ref{norgerm}. We could then realize this by considering the orbifold
map introduced in section \ref{orbmap}.

\medskip

\noindent Patrick Popescu-Pampu \\
Univ. Paris 7 Denis Diderot \hspace{3mm} \\Instit. de
Maths. - UMR CNRS 7586 \\ Equipe "G{\'e}om{\'e}trie et dynamique" \hspace{3mm}\\
Case 7012\\ 
2, place Jussieu  \hspace{3mm}\\ 75251-Paris Cedex 05 \hspace{3mm}\\ FRANCE
\medskip

\noindent ppopescu@math.jussieu.fr

\end{document}